\documentclass[11pt, twoside]{article}

\usepackage{amsmath,amssymb}

\usepackage{bm}

\makeatletter
\@addtoreset{equation}{section}

\makeatother

\begin{document}

\title{\bf Comparison theorems for invariant measures \\

 of random dynamical systems 

}

\author{Tomoki Inoue\\

Graduate School of Science and  Engineering\\
Ehime University\\
Matsuyama 790-8577, Japan\\
E-mail: inoue.tomoki.mz@ehime-u.ac.jp
}

\maketitle

\begin{abstract}
We  study a random dynamical system 
such that one transformation is randomly selected from a family of transformations and then applied on each iteration. 
For such random dynamical systems, 
we consider  estimates of absolutely continuous invariant measures. 
Since the random dynamical systems  are made by  complicated compositions of many deterministic maps
and probability density functions, 
it is  difficult  to estimate the  invariant measures. 
To get rid of this difficulty, we  present fundamental comparison theorems 
 which make easier the estimates of  invariant measures of random maps.
We also demonstrate how to apply the comparison theorems to random maps
with indifferent fixed points and/or  with unbounded derivatives. 
\end{abstract}


\section{Introduction}

One-dimensional piecewise $ C^2$ expanding 
 (Lasota-Yorke type) deterministic  maps have absolutely continuous invariant measures 
whose densities are bounded above ([LY1]) 
since invariant densities are fixed points in a space of functions of bounded variation. 
(In this direction,  there are many studies, for example [BoG] and [LM].)

However,  some  deterministic 
one-dimensional maps with indifferent fixed points have absolutely continuous invariant measures 
whose densities are not bounded above. 
Further, such cases, some estimates of the invariant measures are known.   For example,  see [T1, T2, LSV].
Here, we give an example of such a map.  Let $t \ge 1$ be a constant and let  
$ T: [0,1] \to [0,1]$ be defined by $ T(x) = x + x^t   \mod 1$. 
Then, $ T$ has an absolutely continuous $ \sigma$-finite invariant measure $\mu$ with 
$
\mu([0, \varepsilon]) = \infty 
$
for  $t \ge 2 $, and 
$
\mu([0, \varepsilon]) \approx  \varepsilon^{2-t} 
$
for $ 1 \le t < 2$, 
where $f(\varepsilon)  \approx  g(\varepsilon) $
means that there exist positive constants $c_1$ and $c_2$
such that  $  c_1 g(\varepsilon)  \leq f(\varepsilon)  \leq c_2 g(\varepsilon) $
for any small $\varepsilon > 0 $. 

It is also known that 
one-dimensional piecewise expanding surjective maps satisfying the R\'enyi condition 
have absolutely continuous invariant measures 
whose densities  are bounded away from 0. (For example, see [LM].)
It is also known that 
one-dimensional piecewise expanding maps with some conditions 
have absolutely continuous invariant measures 
whose densities are bounded away from 0 on their supports. 
(For example, see [K] or [BoG].)

Using the  result of [R], we notice  that  some  deterministic 
one-dimensional maps with unbounded derivatives  have absolutely continuous invariant measures. 
However, 
the invariant densities are not bounded away from 0 on the supports, 
which is not necessarily  well-known yet. 
For example, 
$ T: [0,1] \to [0,1]$ defined by
$$
 T (x) =  
 \left\{
\begin{array}{ll} 
\frac 4 3 x \  &   \mbox{ for }  x  \in   \left[0 , \frac 3 4 \right) ,     
\\ 
\left( 4x - 3  \right)^{1/s}    &    \mbox{ for }  x  \in  \left[  \frac 3 4 , 1 \right]
\end{array}  \right.
$$
has  an absolutely continuous invariant probability measure $ \mu$ 
with $ \mu([0, \varepsilon]) \approx \varepsilon^s $  if $1 < s < 4 $, 
which is an offshoot of this research. 
Indeed, this fact follows from Corollary 7.7 in the present paper. 
(If $s> 4$, 1 is an attracting fixed point. 
If $s = 4$, 1 is an indifferent fixed point. 
Hence, we assume $1 < s < 4 $ in the above.) 
The parameter $s$ expresses the  sharpness of $T$ at $ x = 3/4$. 
There are not many general results for maps with unbounded derivatives in this field, 
although statistical properties of some spacial 
intermittent maps  with indifferent fixed points and 
with unbounded derivatives are  studied in  [CHMV]. 
 
Estimates of  absolutely continuous invariant measures of  
deterministic 
one-dimensional maps with indifferent fixed points can be studied using
the induced or the first return maps.  (For example, see [Pi, T1, LSV].)
In this paper,  it is also revealed that 
estimates of  absolutely continuous invariant measures of  
deterministic one-dimensional maps with unbounded derivatives  
can be also studied using
the induced or the first return maps. 

For the random maps  with indifferent fixed points
and/or with unbounded derivatives,
comparing with deterministic maps, 
 it is not easy  to estimate the absolutely continuous invariant measures, 
though we have known in [In4] how to make the first return map of a random map
and how to make an invariant measure of the original random map from that of the first return map. 
In [In4],  we managed to estimate the absolutely continuous invariant measures for some  random maps  with indifferent fixed points.
Since random maps are made by  complicated compositions of many deterministic maps
and probability density functions, 
we need some theorems which make easier the estimates of  absolutely continuous invariant measures of random maps.
For this need, in this paper,  we present fundamental comparison theorems 
 which make easier the estimates of absolutely continuous invariant measures of random maps.

There are many kinds of random dynamical systems  which were studied by many authors. 
For example, 
[BaBD, BaG, GBo, In3, In4, KiLi, MaSU, Mo, Pe] in the references of this paper.
[KiLi] and [MaSU]  are nice surveys on random dynamical systems. 
In [KiLi], there are rich references on random dynamical systems. 

Let explain the random dynamical systems we mainly consider in this paper.
We consider a family of transformations 
$T_t : X \to X \ (t \in W )$ and study a random dynamical system 
such that one  transformation is randomly selected from the family
$\{ T_t : t \in W \} $  and then applied on each iteration.  
In many papers of random maps,  
for example [BaG, BoG, Pe], 
$W$  (the parameter space)   is a finite set.  
 However, in this paper, as in [In3, In4], 
$W$ may even be of cardinality continuum. 
$W$ is not necessary a one-dimensional space. Moreover, we do not assume any 
topological structures on $W$  in our comparison theorems (Theorems 6.2 and 6.5).
Further, in many  papers of random maps, for example 
[Mo, Pe], the selection of $T_t$ is independent of the position  $x \in X$. 
In this paper, 
the selection of $T_t$  may  depend on the position $x \in X$, 
that is, 
$T_t$ is selected according to a probability density function $ p(t, x)$.
We call such a random map a position dependent random map, 
which will be defined precisely in the next section.

Here, we give two examples of random maps 
to which our results in this paper are applicable. 

\vspace{1ex}

\noindent
{\bf Example 1.1. } Let $a_0 $ and $ b_0$ be constants with $1 < a_0 < b_0 $ and let 
$c_0$ be a constant with $0 < c_0 < 1 $ and $ a_0 c_0 \ge 1$. 
Define $ T_{t,s}: [0,1] \to [0,1]$ by 
$$
 T_{t,s} (x) =  
 \left\{
\begin{array}{ll} 
x + \frac 1 t  \left(  \frac t {t-1} \right)^t  x^t  \  &   \mbox{ for }  x  \in   \left[0 , \frac {t-1} t \right) ,     
\\ 
\left( tx - t+ 1  \right)^{1/s}    &    \mbox{ for }  x  \in  \left[  \frac {t-1} t , 1 \right], 
\end{array}  \right.
$$
where $(t, s)$ is randomly selected from 
$ W = \{(t, s) \, | \,  a_0 \le t \le b_0, \,  a_0 c_0  \le s \le c_0 t  \} $ 
accordingly to the  probability density function 
$\displaystyle  p(t, s, x) = 2 / \{c_0 (b_0 - a_0)^2 \}$.

\vspace{2ex}

\noindent
{\bf 
Example 1.2. }
Let $a_0 $, $ b_0$ and $c_0$ be as in Example 1.1.
Define $ T_{t,s}: [0,1] \to [0,1]$ by 
$$
 T_{t,s} (x) =  
 \left\{
\begin{array}{ll} 
\frac t {t-1} x \  &   \mbox{ for }  x  \in   \left[0 , \frac {t-1} t \right) ,     
\\ 
\left( tx - t+ 1  \right)^{1/s}    &    \mbox{ for }  x  \in  \left[  \frac {t-1} t , 1 \right], 
\end{array}  \right.
$$
where $(t, s)$ is randomly selected as in Example 1.1.

\vspace{1ex}

Example 1.1 has an indifferent fixed point 0.
On the other hand, Example 1.2  has no indifferent fixed point.
Both examples have unbounded derivatives at $x =(t-1)/t $.


In these examples, we choose a uniform probability density function which is independent of the position $x$.
If $ p(t, s, x) $ in these examples are replaced with 
$ p(t, s, x) = \frac{ 8 s t (1- x) } { c_0^2 ( b_0^2 - a_0^2)^2 } + \frac {2x } {c_0 (b_0 - a_0)^2 }$,
the random maps really depend on the position $x \in [0,1]$.

The existence of absolutely continuous invariant probability measures of 
position dependent piecewise expanding random maps 
was proved under some conditions  in [In3].  Further, 
the existence of absolutely continuous $\sigma$-finite invariant measures of 
position dependent random maps with indifferent fixed points 
was proved under some conditions  in [In4].

In the present paper,  we assume the existence of an absolutely continuous $\sigma$-finite invariant measure for a random map 
and  give two fundamental  comparison theorems 
 which are  powerful to estimate an absolutely continuous  invariant measure  in a small $\varepsilon$-neighborhood 
 $B(q, \varepsilon) $ of a fixed point  $q$ for  a random map. 
One of the fundamental  comparison theorem is for the lower estimate and the other is for the upper estimate.
These two  theorems are given under general setting and hence 
we can apply these theorems to many kinds of random maps.

In  each comparison theorem,  we compare an absolutely continuous invariant measure
of the original random map with  it of an auxiliary random map.
To use the comparison theorems, we have to select  
 auxiliary random maps whose invariant measures are easily estimated.
Moreover,  we have to select  these to get nice estimates of the invariant measure
of the original random map.

 How do we find the auxiliary random maps?
To answer such a question, 
we will demonstrate how to use our comparison theorems to estimate absolutely continuous invariant measures for random maps 
 with an  indifferent fixed point and/or 
 with unbounded derivatives
whose examples are Examples 1.1 and 1.2. 
Using  our comparison theorems, 
we give some estimates of absolutely continuous invariant measures for these kinds of random maps
(Theorems 7.1 and 7.5).

Remark that 
we are going to  consider the following four other random maps for a given original random map: 

the first return random map, 

the random map with the identity map (This is use in the proof of the comparison theorem for the upper estimate.),

two auxiliary random maps to get nice estimates of invariant measures. 


The paper is organized as follows:
In Section 2, following [In3], we define a position dependent random map $T$
as a Markov process. 
We also study position dependent random maps
with  the identity map.
In Section 3, following [In4], we introduce 
the first return random map of an original random map. 
Further, we summarize how to make an invariant measure of the original random map. 
In Sections 4 and 5, 
in order to feel that the assumptions of our comparison theorems are not too restricted,  
we state on the existence of 
finite and $ \sigma$-finite invariant measures
for one-dimensional random maps.
In section 4, we consider piecewise expanding random maps. 
And in section 5,  we consider random maps 
with indifferent fixed points. 
(The reader who is not interested in such maps  may skip these two sections.)
In Section 6, we present the comparison theorems for invariant measures of random maps. 
In Section 7, we demonstrate how to apply our comparison theorems.
Here, we apply our comparison theorems to random maps with unbounded derivatives and/or  
with indifferent fixed points.
In Section 8, we sate a concluding remark for the possibility of applications  of our comparison theorems.

\vspace{2ex}

\section{Position dependent random maps  }

In the first half of this section, 
following [In3],  
 we define  random maps which may depend on the position in a general setting.  
Further, we  also define invariant measures for  random maps. 
In the latter half of this section, we study position dependent random maps
with  the identity map, which will be used later.

Let 
$( W, {\cal B }, \nu )$ be a $\sigma  $-finite measure space. 
We use $ W $ as a parameter space.
Let 
$ (X, {\cal A },  m) $ be 
a  $\sigma  $-finite measure space.
We use  $X$ as a state space.
Let  
$  T_t : X \to X  $  ($  t \in W$) be 
a nonsingular transformation,
which means that 
$  m(T_t^{-1}D) = 0 $ if $ m(D) = 0 $ for any $ D \in {\cal A } $.
Assume that
$  T_t(x) $ is a measurable  function of $t \in W$ and $x \in X $.

Let $ p : W \times X \to [0, \infty )  $
 be a measurable function  
which is a probability density function of $t \in W$ for each $x \in X $,
that is,     
$  \int_W p(t,x)  \nu (dt) = 1$ for  each $ x \in X  $. 

We define the random map 
$  T = \{ T_t; p(t, x)  : t \in W  \} $
as a Markov process  
with the following transition function:
$$
{\bf P}(x,D) := \int_W  p(t,x)  \, 1_D(T_t (x)) \,  \nu(dt)  \   \  \mbox{for any }   x \in  X \  \mbox{and  for any }  D \in {\cal A },
$$
where  
$ 1_D$ is the indicator function of  $D$.
The transition  function   {\bf P}
induces an operator ${\bf P_{\ast}} $ 
on measures on $X$ defined by 
\begin{eqnarray*}
{\bf P_{\ast}} \mu(D) 
:=  \int_X {\bf P}(x,D) \, \mu(dx) 
= \int_X  \int_W  p(t,x)   \, 1_D(T_t (x)) \,  \nu(dt ) \,  \mu(dx) 
\end{eqnarray*} 
for any measure $\mu$ on $X$ and any $ D \in {\cal A }  $.

If $ {\bf P}_{\ast} \mu = \mu $, $ \mu $
is called an {\em  invariant measure} for
the random map $ T =\{ T_t; p(t,x) : t \in W  \} $. 
If $W$ consists of only one element, $ {\bf P}_{\ast} \mu = \mu $ is 
the usual definition of an invariant measures for a deterministic map.

If $ \mu $ has a density $g$,  then $   {\bf P}_{\ast} \mu $ also has a density, 
which we denote $ {\cal  L}_T g $. 
In other words, $ {\cal  L}_T : L^1(m) \to L^1(m)$ is the operator satisfying 
$$
\int_D  {\cal  L}_T g(x) \, m(dx)   = \int_X  \int_W  p(t,x)   \, 1_D(T_t (x)) \,  \nu(dt ) \,  g(x) \, m(dx) 
$$
for any $ D \in {\cal A }  $. 
We call $ {\cal  L}_T  $ the Perron-Frobenius operator corresponding to 
the random map $ T $. 
If  $W$ consists of only one element, $ {\cal  L}_T  $ is the Perron-Frobenius operator corresponding to 
the deterministic map $ T $.
If you would like to know more the Perron-Frobenius operator 
for a deterministic map, 
see [BoG]  or [LM]. 

Let $  {\cal  L}_{T_t} $ be the Perron-Frobenius operator corresponding to  $ T_t $. 
Then,  we have 
\begin{equation}
{\cal  L}_T g(x)   =   \int_W   {\cal  L}_{T_t} (p(t,x)  g(x))   \,  \nu(dt ).  
\end{equation}

\vspace{2ex}

Now, we prepare  a new random map which we will use 
in the proof of our main result. 
Let $ \hat p : X \to (0,1]$ be a measurable function. 
For the random map $T$ and the probability $ \hat p$, we define  
the new random map 
$ \Upsilon_T= \{T, id; \hat p(x) , 1- \hat p(x) \}$
as a Markov process with the following transition function:
$$
{\bf P}_{\Upsilon_T} (x,D) :=  \hat{p}(x) \int_W p(t,x) \, 1_D(T_t (x)) \,  \nu(dt)  +  (1-   \hat{p} (x) ) 1_D(x)  
$$
for any $ x \in X$ and for any $ D \in  {\cal A}$. 
This is a  random map with the identity map, that is,
$T$ is selected with probability $\hat p (x)$ and $id$ is selected 
with probability $1- \hat p (x)$.

Let $  {\cal  L}_{\Upsilon_T} $ be the Perron-Frobenius operator corresponding to  $ \Upsilon_T$. 
Then, as (2.1),  we have 
\begin{equation}
{\cal  L}_{\Upsilon_T} g (x) 
= {\cal  L}_T \left(  \hat{p}(x) g(x) \right) + \left(1-   \hat{p} (x) \right) g(x).
\end{equation}

For the random map $ \Upsilon_T= \{T, id; \hat p(x) , 1- \hat p(x) \}$, 
we have the following lemma, which will be used later.

\vspace{2ex}

{\bf 
\noindent 
Lemma 2.1.}
{\em 
If the random map $T$ has a $\sigma$-finite invariant measure $\mu$ defined as 
$
\mu (D) = \int_D h(x) \mu (dx) 
$, 
then the random map 
$ \Upsilon_T= \{T, id; \hat p(x) , 1- \hat p(x) \}$ defined above has 
a $\sigma$-finite invariant measure $\hat \mu$ defined as 
$
\hat \mu (D) = \int_D h(x)/\hat p(x)  \mu (dx) 
$. }

\vspace{1ex}

\noindent 
{\it Proof.}
Let $ {\cal L}_{ T} $ and  ${\cal L}_{ \Upsilon_T}$ be the Perron-Frobenius operator 
corresponding to $ T$ and $ \Upsilon_T$,  respectively. 
Then, by (2.2)  we have 
$$
 {\cal L}_{ \Upsilon_T}   \frac { h} { \hat{p} }  
 =  {\cal L}_{ T} \left( \hat{p}  \cdot  \frac { h} { \hat{p} }  \right)
  + (  1-    \hat{p} )  \frac { h} { \hat{p} } 
 =
   {\cal L}_{ T} h +  \frac { h} { \hat{p} } - h 
 = 
  \frac { h} { \hat{p} }. 
$$
This equality implies the lemma.

\vspace{3ex}

\section{First return random maps and invariant measures}

Following [In4], we define 
the first return random map of an original random map. 
Further, we summarize how to make an invariant measure of the original random map 
from that of the first return random  map.

Now, 
let $  T = \{ T_t; p(t, x)  : t \in W  \} $ be a random map on $X$ as in the previous section,  
and 
let $A \subset X$ be a measurable set with $ m(A) > 0$.
If the first return random map on  $A$ of $T$ is well defined,
almost every  $x \in A $ must come back to $A$ by some random iterations.  Thus, 
we assume that 
such a probability is 1, that is,
\begin{equation}
 \sum_{n = 1}^{\infty}  {\hat p}_n (x)  = 1   \quad \mbox{for $m$-a.e.}  x \in A,  
\end{equation}
where 
$$ {\hat p}_1 (x) =  \int_W  p(t,x) \ 1_A  (T_t (x)) \ \nu(dt ),  
$$
$$ {\hat p}_2 (x)  =  \int_{W^2}  p(t_1,x) \   1_{X \setminus A}  (T_{t _1} (x)) \cdot 
p(t_2, T_{t _1} (x) )  \ 1_A  (T_{t _2}  (T_{t _1} (x)) )\ \nu(dt_1 ) \nu(dt_2 ), 
$$
$$
\cdots
$$
\begin{eqnarray*}
  {\hat p}_n (x)
=& &  \hspace{-3ex} \int_{W^n}  p(t_1,x)   \   1_{X \setminus A}  (T_{t _1} (x))  \cdot 
 p(t_2, T_{t _1} (x) )   \   1_{X \setminus A}  (T_{t _2}  (T_{t _1} (x)) )   \\
&  & 
\cdots  \  p(t_{n-1},    T_{t _{n-2}}  \circ  \cdots   \circ   T_{t _1} (x) )   
 \ 1_{X \setminus A}  (T_{t _{n-1}}  \circ  \cdots   \circ   T_{t _1} (x) )   \cdot \\
 &  &
 \  p(t_n,    T_{t _{n-1}}  \circ  \cdots   \circ   T_{t _1} (x) )  
  \ 1_A   (T_{t _{n}}  \circ  \cdots   \circ   T_{t _1} (x) )
 \ \nu(dt_1 ) \nu(dt_2 )  \cdots  \nu(dt_n ).
\end{eqnarray*}
Under this assumption, 
the first return random map $R$ of $T$ on $A$ 
is defined in the following way.

Put 
$$
\widetilde{W} : = W^{{\mathbb N}} 
= \{ (u_1,  u_2,  \ldots )  :   \  u_k \in W,  \  k \in {\mathbb N} \},
$$
which is the parameter space corresponding to the first return random map.
Fix $ x \in A $ arbitrary and put 
\begin{eqnarray*}
\widetilde{W}_1  : =& & \hspace{-3ex}
\{  (u_1,  u_2,  \ldots ) \in \widetilde{W}:   T_{u _1} (x) \in A   \}, 
\\
\widetilde{W}_2  : =& & \hspace{-3ex}
\{  (u_1,  u_2,  \ldots ) \in \widetilde{W}:   T_{u _1} (x) \in X \setminus A,  \ 
     T_{u _2} (T_{u _1} (x))  \in A \}, 
\\
\widetilde{W}_n : =& &  \hspace{-3ex} 
\{  (u_1,  u_2,  \ldots ) \in \widetilde{W} :  T_{u _1} (x) \in X \setminus A,  \ldots ,   \\
&  &  \  
 T_{u _{n-1}}  \circ   \cdots  \circ  T_{u _2}  \circ  T_{u _1} (x)  \in X \setminus A, \
 T_{u _n}  \circ  \cdots  \circ  T_{u _2}  \circ T_{u _1} (x)   \in A \} \ 
    \mbox{ for }  \   n \ge 3. 
\end{eqnarray*}
For  any $ \tilde{u}  =  (u_1,  u_2,  \ldots ) \in \widetilde{W}$,    
we define  $ r_{\tilde{u}} (x) $ and $ Pr(\tilde{u}, x ) $ as follows:

If  $  \tilde{u}  \in \widetilde{W}_1  $,  then $ r_{\tilde{u}} (x) :=  T_{u _1} (x)$
and $ Pr(\tilde{u}, x ) :=  p(u_1, x)$.

If  $  \tilde{u}  \in \widetilde{W}_2 $,
then 
$ r_{\tilde{u}} (x) :=  T_{u _2} (T_{u _1} (x))$ and
$ Pr(\tilde{u}, x ) :=  p(u_1, x) \cdot p(u_2, T_{u _1} (x) ) $. 

If $  \tilde{u}  \in \widetilde{W}_n  $,  
then  $   r_{\tilde{u}} (x) :=   T_{u _{n}}  \circ   \cdots \circ   T_{u _2} \circ  T_{u _1} (x) $  and 
$$
 Pr(\tilde{u}, x ) :=  p(u_1, x) \cdot p(u_2, T_{u _1} (x) )    \cdots     
p(u_n,   T_{u _{n-1}}  \circ   \cdots   T_{u _2}  \circ  T_{u _1} (x) ) 
\    \mbox{ for }  \   n \ge 3. 
$$

\vspace{2ex}

Define the measure $ \tilde{\nu} $ as 
$$
 \tilde{\nu} (d \tilde{u}) = \sum_{n=1}^{\infty} 1_{\widetilde{W}_n}  (\tilde{u})
     \nu(du_1) \cdots  \nu(du_n).
$$
Then, by (3.1),   
$  Pr(\tilde{u}, x ) $ is a probability density function of $\tilde{u} $
for $m$-a.e. $x \in A$, that is,  
$$
\int_{\widetilde{W}}   Pr(\tilde{u}, x )  \tilde{\nu}(d\tilde{u})  
=  
\sum_{n=1}^{\infty}  \int_{{\widetilde{W}}_n }  Pr(\tilde{u}, x )  \tilde{\nu}(d\tilde{u})  
= 
\sum_{n = 1}^{\infty}  {\hat p}_n (x)  = 1
  \quad \mbox{ for $m$-a.e.} \, x \in A.
$$

So, 
the first return random map $  R = \{  r_{\tilde{u}} \  ; Pr(\tilde{u}, x ) :   \tilde{u} \in \widetilde{W} \} $ 
of the random map $T$ on $A$  
is defined as a Markov process  
with the following transition function:
$$
{\bf P}_R  (x,D) := \int_{\widetilde{W}}   Pr(\tilde{u}, x )  \ 1_D(r_{\tilde{u}} (x)) \  \tilde{\nu}(d\tilde{u})
$$
for any $x \in A $ and for any $D \in  {\cal A}   $  with  $D \subset A$.
In the same way as  the transition function $ {\bf P} $ 
induces an operator ${\bf P}_{\ast} $, 
the transition  function   ${\bf P}_R$  
induces an operator ${\bf P}_{R_{ \ast}} $ 
on measures on $A$ defined by 
\begin{eqnarray*}
{\bf P}_{R_{ \ast}} \mu_A(D) 
: &= & \int_A {\bf P}_R  (x,D) \mu_A (dx) \\
&=  & \int_A  \int_{\widetilde{W}}  
             Pr(\tilde{u}, x )  \ 1_D(r_{\tilde{u}} (x)) \  \tilde{\nu}(d\tilde{u}) \  \mu_A (dx)
\end{eqnarray*} 
for any measure $ \mu_A $ on $A$ and for any $ D \in  {\cal A}  $ with $ D \subset A $. 
Thus, by the definition of the invariant measure of a random map, 
 $ \mu_A $ is an invariant measure for the first return random map $ R $ 
if  $  {\bf P}_{R_{ \ast}} \mu_A = \mu_A$.

\vspace{1ex}

To state how to make an invariant measure of the original random map, 
we define  two operators  ${\cal U}_T$ and $ \widetilde{{\cal U}}_T $  : $ L^{\infty} (m)  \to  L^{\infty} (m)$ 
by 
$$
{\cal U}_T f(x) = \int_W  p(t, x)  \  f(T_t(x))  \   \nu (dt)  \   \ 
\mbox{and}  \   \  
\widetilde{{\cal U}}_T  f =  {\cal U} (1_{X \setminus A }  \cdot  f )    \   \ 
  \mbox{for}   \  \   f   \in L^{\infty} (m).
$$

\vspace{1ex}

If $T_0 = T_t $ for any $ t \in W $ ($T $ is deterministic), 
$  {\cal U}_T f(x) = f(T_0 (x))$, that is,  
$  {\cal U}_T $ is  the Koopman operator corresponding to $ T = T_0 $. 
By this reason,  for a random map $T$, $ {\cal U}_T$ defined above  should be called 
the random Koopman operator. 

\vspace{1ex}

We summarize the relations  
among some symbols appeared in the present and the previous sections:

\vspace{2ex}

\noindent 
 {\bf  Proposition 3.1. } (Proposition 4.1 in [In4])
 {\em Assume that the first return random map on $A$ of $T$ is well defined, 
that is, (3.1)  is satisfied. 
  Let ${\bf P}_R $, ${\cal U}_T$ and $ \widetilde{{\cal U}}_T $ be defined above.
Let 
$  {\bf P} $ and $  {\bf P}_{\ast}  $ be as in Section 2. 
Then,

\begin{itemize}
\item[{\rm  (i)} ]
$ \displaystyle 
{\cal U}_T 1_D (x) =   {\bf P}(x,D)
$ 
for any $  D \in {\cal A } $. 

\item[{\rm  (ii)} ]
$ \displaystyle 
 {\bf P}_{\ast} \mu = \mu  $ \   if and only if   \  $ \int_X  {\cal U}_T f (x) \   \mu (dx)   =  \int_X  f(x) \  \mu (dx)  \ 
 \mbox{for any}  \   f  \in  L^{\infty} (m).
$

\item[ { \rm (iii)} ]
$\displaystyle 
{\hat p}_{n} (x) = { \widetilde{{\cal U}}}_T^{n-1}  \left(  {\cal U} 1_A   \right) (x)
$ 
for $ n \in {\mathbb N} $, 
where 
$\displaystyle    { \widetilde{{\cal U}}}_T^{0}  \left(  {\cal U}_T 1_A   \right) 
=  {\cal U}_T 1_A   $.

\item[  {\rm (iv)}  ]
$ \displaystyle 
{\bf P}_R  (x,D) =  \sum_{n = 0}^{\infty}  { \widetilde{\cal U}}_T^n  \left(  {\cal U}_T 1_D   \right) (x)
$   
for any $  D \in {\cal A } $ with $D \subset A  $. 

\end{itemize}

}

\vspace{2ex}

From [In4], we know the following two theorems,
which make clear the relation between an invariant measure of the original random map 
and that of the first return random map.

\vspace{3ex}

\noindent 
{\bf Theorem 3.2.}
{\em 
Let  $T$ be a random map as in Section 2.
Assume that  $T$ has a $  \sigma $-finite  invariant measure $ \mu  $ 
and that the  first return random map $R$  of $T$ on $A$ is well defined.
Then,  $ \mu|_A $ is an invariant measure for  $R$, 
where $ \mu|_A  $ is the restriction of $ \mu$ to $A$.  
}

\vspace{3ex}

\noindent 
{\bf Theorem 3.3.} 
{\em 
Let  $T$  be a random map as in Section 2.
Assume that the  first return random map $R$  of $T$ on $A$ is well defined.  Further,  
assume that $R$ has an invariant probability measure $ \mu_A $.
Define the measure $ \mu $ as 
\begin{equation}
\mu (D)  
=   \sum_{n = 0}^{\infty}    \int_A  { \widetilde{{\cal U}}}_T^n  \left(  {\cal U}_T 1_D   \right) (x) \  \mu_A (dx) 
\quad \mbox{for any}  \  D \in {\cal A}.   
\end{equation}
Then,  $\mu$ is a $  \sigma$-finite invariant measure for $T$.
}

\vspace{2ex}

\section{Invariant measures of piecewise expanding random maps}

In sections 2 and 3, the settings of random maps are abstract. 
From now on, we consider one-dimensional random maps. 
In this section, we summarize the existence of absolutely continuous invariant measures 
for piecewise expanding random maps.

Set $ X=  [0,1]  $ and set $m$  the 
Lebesgue measure.  Then, 
$ T_t  $ 
is a map from $ [0,1] $ into itself
for each $t \in W$,
$ T_t(x) $ is a measurable  function of $t \in W$ and $x \in [0,1] $,
and 
 $ p : W \times [0,1] \to [0, \infty )  $
is a measurable function  
such that      
$  \int_W p(t,x)  \nu (dt) = 1$ for each $ x \in [0,1]  $. 

Let 
$ \Lambda $ be  a countable or finite set
and let 
 $\Lambda_t \subseteq \Lambda $ 
for each $ t \in W $.
We use $ \Lambda $ as a set of indices of subintervals of $[0,1] $.
For each $ t \in W $, 
we assume that
$\{ I_{t,i}  \}_{i \in \Lambda_t}$
is a  family of closed intervals
 such that  
${\rm int} (I_{t,i}) \cap {\rm int} (I_{t,j})  
= \emptyset
$ ($ i \ne j $ )  and 
$  m(  [0,1] \setminus \cup_{i \in \Lambda_t}   I_{t,i}) =0  $, 
where $ {\rm int}(I) $ is the interior of an interval $I$.

To avoid  using the subscript $t$ on $ \Lambda_t $,
 we set  $ I_{t,i} = \emptyset $ for $ i \in  \Lambda \setminus \Lambda_t $.
For  convenience, we consider the empty set  as a closed interval.

For a random map 
$ T = \{ T_t ; p(t,x);  \{ I_{t,i}  \}_{i \in \Lambda}  : t \in W \} $,
we make  two basic assumptions 
and  introduce piecewise $C^1$ random maps.

The first assumption is :
%
%
\begin{enumerate}
\item [(A1)]   
The restriction to $ {\rm int}  (I_{t,i}) $ of $T_t $ is a 
$ C^1$ and monotone function for each $ i \in  \Lambda $ and $ t \in W $.
\end{enumerate}


Let 
$ T_{t,i}$ be  the restriction of $ T_t  $ to ${\rm int}  ( I_{t,i}  ) $ 
 for each $ t \in W $ and $ i \in  \Lambda  $.
Put 
$$ 
 \phi_{t,i}(x) : = 
\left\{ \displaystyle
\begin{array}{ll} \displaystyle
T_{t,i}^{-1}(x),   &    \quad  x \in T_{t,i} ( {\rm int} (I_{t,i}) ),  \\
0, &  \quad  x \in [0, 1] \setminus T_{t,i} ( {\rm int} (I_{t,i}) ) 
\end{array}
 \right.   
$$
for each $ t \in W $ and $ i \in  \Lambda  $.
We note that $ \phi_{t,i}(x) = 0 $ 
if $ i  \in  \Lambda \setminus \Lambda_t $.

The second assumption is :

\begin{enumerate}
\item[(A2)]  
For each $x \in X$ and $i \in \Lambda $,
$ w_{x,i}(t) := \phi_{t,i}(x) $
is a measurable function of $t$.
\end{enumerate}


We call a random map $T = \{ T_t, p(t, x),   \{ I_{t,i}  \}_{i \in \Lambda}  : t \in W  \} $
 satisfying the assumptions  (A1) and  (A2) 
a {\em piecewise $C^1 $ random map}.

From  [In3], we know the following theorem.

\vspace{2ex}

\noindent 
{\bf Theorem 4.1.  }
{\em 
Let $T = \{ \tau_t;  p(t, x) ;  \{ I_{t,i}  \}_{i \in \Lambda}  : t \in W  \} $ be 
a piecewise $C^1 $ random map. 
Put 
\begin{eqnarray}
g(t,x) = 
\left\{ \displaystyle
\begin{array}{ll} \displaystyle
\frac{p(t,x)}{| \tau_t^{\prime } (x) |},   &    \quad  x \in \cup_i {\rm int} (I_{t,i}), \\
0,   &  \quad  x \in [0, 1] \setminus \cup_i {\rm int} (I_{t,i}) 
\end{array}
 \right.   
\end{eqnarray}
for $ t \in W $.
Suppose  the following conditions:
\begin{itemize}
\item[{\rm (a)}] 
$  \sup_{x \in [0,1]}  \int_{ W}  g(t,x) \  \nu (dt) 
 < 1 $.

\item[{\rm (b)}] 
There exists a constant  $M$ 
such that 
$ \bigvee_{[0,1]} g(t, \cdot )  < M $ 
for a.s. $t \in W$, where 
$  \bigvee_{[0,1]} g(t, \cdot ) $ stands for the total variation of $ g(t, \cdot ) $ 
on $[0, 1] $.
\end{itemize}
Then,
 the random map $T$  has an absolutely continuous invariant probability measure.
Moreover, the density of the invariant measure is of bounded variation.

}

\vspace{2ex}

\section{$ \sigma $-finite invariant measures of random maps with indifferent fixed points}  

In [In4], 
we obtained  the following result  on the existence of 
an absolutely continuous $ \sigma $-finite invariant measures 
for a piecewise $C^1 $ random map with indifferent fixed points.

\vspace{2ex}

\noindent
{\bf Theorem 5.1.}
{\em 
Let 
$ T = \{ T_t ; p(t,x);  \{ I_{t,i}  \}_{i \in \Lambda}  : t \in W \} $ be a piecewise $ C^1 $ random map
defined in Section 4.
For $ t \in W $ and $ x \in [0,1] $, 
put $ g(t,x) $ as in (4.1).
Let $  x_1, x_2, \ldots, x_{n_0} $ be some common fixed points of all $ T_t   \  (t \in W) $. 
For each $ k  \in \{ 1, 2, \ldots, n_0  \}$,   let $ U_k$ be 
a small one sided or two sided neighborhood of  $ x_k $. 
Put $ A = [0,1] \setminus  \cup_{k = 1}^{n_0} U_k$.
Suppose  the following conditions:
\begin{itemize}
\item[{\rm (i)}]  $  \sup_{x \in A}  \int_W  g(t, x)  \nu (dt) <   1$. 

\item[{\rm (ii)}]    
For each $ t \in W, $
$ g(t, x) $ is  monotonically increasing as $x  \in  U_k $  goes to $ x_k$.

\item[{\rm (iii)}] There exists a constant  $M$ 
such that $ \bigvee_{[0,1]} g(t, \cdot )  < M $ 
for a.s. $t \in W$,
that is,  
there  exists a $ \nu$-measurable set 
$W_0 \subset W   $
such that $ \int_{W_0} p(t,x) \nu (dt) = 1 $ and that 
$ \bigvee_{[0,1]} g(t, \cdot )  < M $ for  all $t \in W_0 $.

\item[{\rm (iv)}]  
For each  $t \in W,$
there exists a couple of  bounded constants  $  m_{t, k} > 0  $ and  $ d_{t,k} > 1$ such that 
 $$
  | T_t(x) - x_k  |  =  |  x - x_k | +   m_{t, k}  |  x - x_k |^{d_{t,k}} +  o(|  x - x_k |^{d_{t,k}} ) 
  \ 
 \mbox{  in }  \ U_k.
 $$

\end{itemize}
Then, there exists a $ \sigma$-finite $T$-invariant measure $ \mu $
which is absolutely continuous with respect to 
the Lebesgue measure. Further, $  \mu ( A)  < \infty $, 
and the density of $\mu  $  is of bounded variation  on $  A$. 

}

\vspace{3ex}

In this theorem, 
$ x_i $ may even equal to $ x_j $ for some $  j \ne i $.
In this case, 
one of $ U_i $ and $ U_j $ should be  a small  left hand side neighborhood of $ x_i = x_j$
and the other should be  a  small right hand side neighborhood of it.  

We remark that 
$  \int_W  g(t, x_k)  \nu (dt) = 1$ for each $ k  \in \{ 1, 2, \ldots, n_0  \}$, 
which means $ x_k$ is an indifferent fixed point for each $k$.

\vspace{3ex}


\section{Comparison theorems for invariant measures}

Let 
$ T = \{ T_t ; p(t, x); \{ I_{t,i}  \}_{i \in \Lambda}  : t \in W \} $ be a one-dimensional 
random map. 
We often assume the existence of an absolutely continuous invariant measure of $T$, 
which is not too restricted as we know the facts of the previous sections.

For simplicity, we only consider 
a small neighborhood of 
the common fixed point 0 of $ T_t  \  (t \in W)$,  
since we can easily understand that similar results hold for a small neighborhood of any common fixed point
in a similar situation.

Let $ \varepsilon_0 > 0$ be a sufficiently small constant and $B(0) = [0,  \varepsilon_0]$. 
We suppose  the following assumtion (6A).

\begin{itemize}

\item[(6A)]

\begin{itemize}
\item [(i)]
$T_t(0) = 0$ for each  $ t \in W $.  

\item [(ii)]
 $(0,  \varepsilon ) \subsetneq  T_t  (0,  \varepsilon ) $ for any $ 0 < \varepsilon  < c $, 
 where $c$ is a constant with $0 < c < 1 $.

\item [(iii)]
The random map $ T = \{ T_t , p(t, x);  t \in W \} $ 
has 
an absolutely continuous $\sigma$-finite invariant measure $  \mu $ with $\mu([c,1]) =1 $. 
Let $h$ be a density of $  \mu $.
Then, 
$$
\gamma_1 \le h(x) \le \gamma_2 
\quad \mbox{on $\displaystyle   \bigcup_{ t \in W} T_{t}|_{[c,1]}^{-1}(B(0)) $}, 
$$
where $\gamma_1  $ and $ \gamma_2 $ are positive constants. 

\end{itemize}

\end{itemize}

\vspace{2ex}

\noindent 
{\bf Remark 6.1.}
To clarify the estimate of the invariant measure,  we assume  $\mu([c,1]) =1 $ in (6A). 
If there is an indifferent fixed point in $ [c,1]$, it may happen that  $\mu([c,1]) = \infty $.  
In such a case, instead of the assumption $\mu([c,1]) =1 $,  
we may assume $\mu(A) =1 $, 
where  $A$ is a suitable measurable set
including the set $  \bigcup_{t \in W} T_t|_{[c,1]}^{-1}(B(0))  $.

\vspace{2ex}

{\bf
\noindent  6.1. 
 Lower estimate. 
 }

Here, we consider the lower estimate 
of $ \mu ([0, \varepsilon ])$
under the following assumption (6.1A). 

\begin{itemize}

\item[(6.1A)]

\begin{itemize}

\item[(i)] 
There exists a function 
$ \tau_1 : [0, c) \to [0,1 )$ which is  monotone  and  satisfies the following condition:

there exists a constant $c_{\ast} $  with $0 <  c_{\ast} \le  \min \{ c , \varepsilon_0 \}$ such that 
 $ T_t  (0,  \varepsilon ) \subset \tau_1 (0,  \varepsilon ) $ 
for any $ \varepsilon  \in (0, c_{\ast})$ and for a.s. $t \in W $,
which means 
$$\int_{W_1(\tau_1)}  p(t,x) \nu(dt) = 1  \quad \mbox{for any $ x \in  [0, c)$},
$$
where 
$ W_1 (\tau_1)
= \{ t \in W \  | \  T_t  (0,  \varepsilon ) \subset \tau_1 (0,  \varepsilon ) \ \mbox{for any $ \varepsilon  \in (0, c_{\ast})$}
\}
$.

\item[(ii)] 
For each $t \in W$,
define $ \bar{\tau}_t : [0, 1]  \to [0,1]$ 
by 
$$
\bar{\tau}_t (x) := 
\left\{
\begin{array}{ll}  
\tau_1 (x)  &    \quad \mbox{for} \ x \in [0, c),  \\   
T_t(x) &  \quad  \mbox{for} \ x \in [ c, 1].
\end{array}
 \right.
$$

The random map 
$ \bar{\tau} = \{ \bar {\tau}_t ;  p(t,x) : t \in W \}$
has an absolutely continuous invariant measure $ \bar {\mu}$ with $\bar \mu([c,1]) =1 $.
Let $\bar h$ be a density of $ \bar {\mu}$. Then $\bar h$ satisfies 
$  \bar{h} (x) \le  \bar \gamma_2 < \infty$ on $  \bigcup_{t \in W} T_t|_{[c,1]}^{-1}(B(0)) $, 
where $ \bar \gamma_2$ is a constant.

\end{itemize}

\end{itemize}

We will show the following theorem, which would be  natural. 

\vspace{2ex}
 
 \noindent 
{\bf Theorem 6.2.}  {\em
Let $ T = \{ T_t ; p(t, x); \{ I_{t,i}  \}_{i \in \Lambda}  : t \in W \} $ be a one-dimensional 
random map which satisfies (6A) and (6.1A).
Then, 
$$ \displaystyle 
  \frac{\gamma_1 }{ \bar \gamma_2 } \bar{\mu}( [0, \varepsilon]  ) \leq  
 \mu([0, \varepsilon ]) 
\ \mbox {for any }  \varepsilon  \in (0, c_{\ast}).
$$

}

\vspace{3ex}

\noindent 
{\bf Remark 6.3.} \ 
The random map $  \bar{\tau}$ is an auxiliary random map 
to  get a lower estimate of $  \mu([0, \varepsilon ]) $. 
If we choose a good  $ \tau_1$, we can obtain a good lower estimate of $  \mu([0, \varepsilon ]) $. 
For example, if we would like to show $  \mu([0, \varepsilon ]) = \infty $, 
we should find $ \tau_1$ such that the random map $ \bar \tau $ has the invariant measure $ \bar \mu  $
 with $ \bar{\mu}( [0, \varepsilon]  ) = \infty$.

\vspace{2ex}

\noindent 
{\bf Remark 6.4.} \ 
We do not need the lower estimate $ \bar h \ge \bar \gamma_1 >0 $ to prove Theorem 6.2. 
However, we use this for the lower estimate of $\bar{\mu}( [0, \varepsilon]  ) $ 
in many cases of applications. 
(See the proof of Theorem 7.1.)

\vspace{2ex}

\noindent 
{\em  Proof of Theorem 6.2. }
 
Let $ R $ be the first return random map of $T$ on $A = [c, 1]$. 
As we  studied  in Section 3, 
the restriction of $h$   to $A$ is  a stationary density of $R$ and 
$  \mu_A = \mu |_A $ is an $ R $-invariant probability measure.
For $T$,  we define 
two operators  ${\cal U}_T$ and $ \widetilde{{\cal U}} _T$  : 
$\displaystyle  L^{\infty} (m)  \to  L^{\infty} (m)$ 
as in Section 3. 
We also define 
two operators  ${\cal U}_{\bar \tau}$ and $ \widetilde{{\cal U}} _{\bar \tau}$  : 
$\displaystyle  L^{\infty} (m)  \to  L^{\infty} (m)$ 
for the random map $ \bar \tau$.

Since $ h \ge  \gamma_1 $ on $  \bigcup_{t \in W} T_t|_{[c,1]}^{-1}(B(0))  $, 
by Theorem 3.3  we have 
\begin{eqnarray}
&  & \mu ( [ 0, \varepsilon]  )
=  \sum_{n = 0}^{\infty}    \int_A
             { \widetilde{{\cal U}}}_T^n  \left(  {\cal U}_T 1_{  [0, \varepsilon] }   \right) (x) \  \mu_A (dx) 
= \sum_{n = 0}^{\infty}    \int_A
             { \widetilde{{\cal U}}}_T^n  
             \left(  {\cal U}_T 1_{  [0, \varepsilon] }   \right) (x) \cdot  h(x)  \  m |_A (dx)       
              \nonumber              \\
&  & \ge    \gamma_1  \   \sum_{n = 0}^{\infty}    \int_{ A }
             { \widetilde{{\cal U}}}_T^n  \left(  {\cal U}_T 1_{  [0, \varepsilon]}   \right) (x) \  m |_A (dx). 
\end{eqnarray}
For any small $  \varepsilon  \in (0, c_{\ast})$ it follows from  the assumption  that 
\begin{eqnarray*}
&  &   \int_{ A }
             { \widetilde{{\cal U}}}_T  \left(  {\cal U}_T 1_{  [0 , \varepsilon]  }   \right) (x) \  m|_A (dx)  \\
&  = &
 \int_{W^2}  \left( 
         \int_{A}   
        p(t_0, x)  \, 1_{[ 0,  {T_{t_1} |_{[0, c) }}^{-1} (\varepsilon)  ] } 
                         ( T_{t_0} (x)) \  p(t_1, T_{t_0} (x))   \, m(dx)  
         \right) \nu (dt_0) \nu (dt_1)    
\\ 
&  \ge &
 \int_{W^2}  \left( 
         \int_{A}   
        p(t_0, x)  \, 1_{ [ 0, \tau_1^{-1} (\varepsilon) ] } 
                         ( T_{t_0} (x)) \  p(t_1, T_{t_0} (x))   \, m(dx)  
         \right) \nu (dt_0) \nu (dt_1)  
\\
&  = &
    \int_{W}  \left( 
         \int_{A}   
        p(t_0, x)  \, 1_{ [ 0,  \tau_1^{-1} (\varepsilon) ] } 
                         ( T_{t_0} (x))    \, m(dx)  
         \right) \nu (dt_0)     
\\
 &  = &
   \int_{ A }
             { \widetilde{{\cal U}}_{\bar \tau} }  \left(  {\cal U}_{\bar \tau}   1_{  [0, \varepsilon]  }   \right) (x) \  m|_A (dx).        
\end{eqnarray*}
Similarly we obtain
$$
 \int_{A}
             { \widetilde{{\cal U}}}_T^n  \left(  {\cal U}_T 1_{  [0, \varepsilon]  }   \right) (x) \  m|_A (dx)
 \ge 
\int_{A}   { \widetilde{{\cal U}}}_{\bar \tau}^n  \left(  {\cal U}_{\bar \tau}  1_{  [0, \varepsilon]  }   \right) (x) \  m|_A (dx),         
$$
and hence we have 
\begin{eqnarray}
 \sum_{n = 0}^{\infty}    \int_{A}
             { \widetilde{{\cal U}}}_T^n  \left(  {\cal U}_T 1_{  [0, \varepsilon]  }   \right) (x) \  m|_A (dx)
 \ge   
 \sum_{n = 0}^{\infty}    \int_A
             { \widetilde{{\cal U}}}_{\bar \tau}^n  \left(  {\cal U}_{\bar \tau}  1_{  [0, \varepsilon]  }   \right) (x) \  m|_A (dx).                  
\end{eqnarray}
Since $ \bar h \le  \bar \gamma_2 $  on $  \bigcup_{t \in W} T_t|_{[c,1]}^{-1}(B(0))  $, we have 
\begin{eqnarray}
 \sum_{n = 0}^{\infty}    \int_A
             { \widetilde{{\cal U}}}_{\bar \tau}^n  \left(  {\cal U}_{\bar \tau}  1_{  [0, \varepsilon]  }  \right) (x) \ m|_A (dx) 
 & \ge  &  \frac 1 {\bar \gamma_2 }    
  \sum_{n = 0}^{\infty}    \int_A
             { \widetilde{{\cal U}}}_{\bar \tau}^n  \left(  {\cal U}_{\bar \tau}  1_{  [0, \varepsilon]  }   \right) (x) \  \bar{\mu}|_A (dx) 
      \nonumber \\
& =  &   \frac 1 {\bar \gamma_2 }     \    \bar{\mu} ( [0, \varepsilon]  ).
\end{eqnarray}

Therefore, by (6.1),  (6.2) and (6.3) we obtain 
$$
 \mu ( [ 0, \varepsilon]  ) \geq \frac {\gamma_1}{\bar \gamma_2 }  \bar{\mu} ( [0, \varepsilon]  )
 \quad \mbox{ for any  $ \varepsilon  \in (0, c_{\ast})$}, 
$$
which shows the theorem.  \hfill $\Box$

\vspace{3ex}

{\bf 
\noindent 
6.2. Upper estimate.}

Next, we consider the upper estimate of $ \mu ([0, \varepsilon ])$
under the following assumption (6.2A).

\begin{itemize}

\item[(6.2A)]

\begin{itemize}

\item[(i)]  There exists a function 
$ \tau_2 : [0, c) \to [0,1 )$ which is  monotone and  satisfies the following condition:

there exists a constant $c_{\ast} $ with $0 <  c_{\ast} \le  \min \{ c , \varepsilon_0 \}$ such that 
$$ W_2 (\tau_2)
= \{ t \in W \, | \, 
T_t  (0,  \varepsilon ) \supset \tau_2 (0,  \varepsilon ) \ \mbox{for any $ \varepsilon  \in (0, c_{\ast})$}
\}
$$
is a measurable set and satisfies 
$$
\hat p (x) :=
\int_{W_2(\tau_2)}  p(t, x) \nu(dt) > 0 
\quad \mbox{for any $ x \in  [0, c) $.}
$$

\item[(ii)] 
For each $t \in W$,
define $ \hat{\tau}_t : [0, 1]  \to [0,1]$ 
by 
$$
\hat{\tau}_t (x) := 
\left\{
\begin{array}{ll}  
\tau_2 (x) &    \quad \mbox{for} \ x \in [0, c),  \\   
T_t(x) &  \quad \mbox{for} \ x \in [ c, 1].
\end{array}
 \right.
$$
The random map 
$ \hat{\tau} = \{ \hat{\tau}_t , p(t,x) ; t \in W \}$
has an absolutely continuous invariant measure $ \hat{\mu}$
with $ \hat \mu([c,1]) =1 $.
Let $\hat h$ be a density of $ \hat{\mu}$. 
Then $\hat h$ satisfies   $ 0< \hat \gamma_1  \le   \hat{h}(x) $
 on $ \bigcup_{t \in W} T_t|_{[c,1]}^{-1}(B(0)) $, 
where $ \hat \gamma_1$  is a constant.

\end{itemize}

\end{itemize}

On the upper estimate of $\mu ([0, \varepsilon ]) $ we are going to show the following theorem.

\vspace{2ex}


{\noindent 
\bf Theorem 6.5.
}
{\em 
Let $ T = \{ T_t ; p(t, x); \{ I_{t,i}  \}_{i \in \Lambda}  : t \in W \} $ be a one-dimensional 
random map which satisfies (6A) and (6.2A). 
Then, 
$$ \displaystyle 
\mu([0, \varepsilon ]) 
\leq   \frac{\gamma_2 }{ \hat \gamma_1  \inf_{x \in [0,c_\ast)} \hat p (x)} \hat{\mu}( [0, \varepsilon]  )
\mbox{ for any  } \varepsilon  \in (0, c_{\ast}).
$$

}

\vspace{3ex}

\noindent 
{\bf Remark 6.6.} \ 
The random map $  \hat{\tau}$ is an auxiliary random map to get an upper estimate of $  \mu([0, \varepsilon ]) $. 
If we choose a good  $ \tau_2$, we can obtain a good upper estimate of $  \mu([0, \varepsilon ]) $. 
In many cases, if we need a good upper estimate, 
 $ \tau_2$ should be given slightly smaller than $ \sup_{t \in W} T_t(x)$
 on $( 0, c_\ast )$. 
If there exists a $ t_0 \in W$ such that 
$ T_t(x) \le T_{t_0} (x)$  and $p(t_0, x) \nu ({t_0}) > 0$ on $( 0, c_\ast )$, 
then we  should  select $\tau_2 = T_{t_0} $ on  $[ 0, c )$.

\vspace{2ex}

\noindent 
{\bf Remark 6.7.} \   Similarly to Remark 6.3, 
we do not need the upper  estimate $ \hat h \le \hat  \gamma_2  < \infty $ to prove Theorem 6.5. 
However, we  use this for the upper estimate of $\hat{\mu}( [0, \varepsilon]  ) $ in many cases. 
(See the proof of Theorem 7.1.)

\vspace{2ex}

To prove  Theorem 6.5 we need a  random map with the identity map. 
To make this new random map,  we 
define the function $ \hat{p} $ on $[0, 1] $ by 
$$
 \hat {p} (x) = 
 \left\{
\begin{array}{ll} 
\int_{W_2(\tau_2)}  p(t, x) \nu(dt)   \  &   \mbox{ for }  x  \in  [0 , c],     
\\ 
1    &   \mbox{ otherwise. } 
\end{array}  \right.
$$
Using $  \hat{p} $, 
we define  the random map $ \Upsilon_{  \hat{\tau}} = \{   \hat{\tau} , id ;  \hat{p}(x),  1-   \hat {p} (x)   \}  $ 
as a Markov process with the following transition function:
$$
{\bf P_{ \Upsilon_{  \hat{\tau}} }  }(x,D) :=  \hat{p}(x) \int_W p(t,x)1_D(\hat{\tau}_t (x)) \nu(dt)  +  (1-   \hat{p} (x) ) 1_D(x)  
$$
for any $x \in [0,1] $ and  for any $  D \in  {\cal A}$.

\vspace{1ex}

The following lemma is an immediate consequence of Lemma 2.1.

\vspace{2ex}

{\bf 
\noindent Lemma 6.8.}
{\em 
Let 
the measure $  \hat{  \hat{\mu}} $ be defined as 
$   \hat{  \hat{\mu}} (D) = \int_D { \hat{ h } (x)}/ { \hat{p}(x) } m(dx) $ $\mbox{for any  }  D \in  {\cal A}$.
Then, $ \hat{  \hat{\mu}} $
 is an absolutely continuous invariant measure
for the random map $  \Upsilon_{  \hat{\tau}} = \{    \hat{\tau} , id ;   \hat{p} (x) ,  1-    \hat{p} (x)   \}  $. 
}

\vspace{2ex}

In the following lemma and the proof of the theorem, 
for simplicity, we omit  $  \hat{\tau}$ of $  \Upsilon_{  \hat{\tau}} $.

\vspace{2ex}

{\bf
\noindent  Lemma 6.9.} 
{\em 
For the random map $  \Upsilon$ mentioned above, 
we define 
two operators  ${\cal U}_{  \Upsilon}$ and $ \widetilde{{\cal U}}_{  \Upsilon}  $ $  : $
$\displaystyle  L^{\infty} (m)  \to  L^{\infty} (m)$ 
as in Section 3.
Then,  for any $x \in [c,1]$ the following two equalities hold:

{\rm (i)} $\   \displaystyle 
{ \widetilde{{\cal U}}}_{  \Upsilon } f (x)
= 
   \int_W p(t_0, x) 1_{   [0 ,  c)   } (\hat \tau_{t_0} (x)) \,
f   (\hat \tau_{t_0} (x)) \nu(dt_0)
$

\vspace{1ex}

{\rm (ii)}  \ 
$ \displaystyle 
 { \widetilde{{\cal U}}}_{  \Upsilon }   \left(  {\cal U}_{  \Upsilon} 
   1_{   [0 , \varepsilon]   }   \right) (x)
$

\qquad $ \displaystyle 
=  \int_W p(t_0, x)
 \left\{
  \hat p (T_{t_0} (x) )  1_{   [0 , \tau_2^{-1} (\varepsilon) ]   }   (T_{t_0}(x))
   + \left( 1 - \hat p (T_{t_0} (x)) \right) 1_{   [0 ,   \varepsilon]   } (T_{t_0} (x))
 \right\} \nu(dt_0).
$

}

\vspace{2ex}

\noindent 
{\it Proof.}

(i)
By the definition of $ { \widetilde{{\cal U}}}_{  \Upsilon }$, we have
\begin{eqnarray*}
&  &  \hspace{-5ex}
 { \widetilde{{\cal U}}}_{  \Upsilon }  f (x) \\
&=&
\hat p (x)  \int_W p(t_0 , x) 1_{   [0 ,  c)   } (\hat \tau_{t_0} (x)) \,
f   (\hat \tau_{t_0} (x)) \nu(dt_0)
  + \left( 1 - \hat p (x) \right)
  1_{   [0 ,  c)   } (x) \, 
f  (x) \\
&=&
   \int_W p(t_0, x) 1_{   [0 ,  c)   } (\hat \tau_{t_0} (x)) \,
f   (\hat \tau_{t_0} (x)) \nu(dt_0)
 \quad \mbox{ for any $x \in [c,1]$. }
\end{eqnarray*}

(ii)
Since we have 
$$
  {\cal U}_{  \Upsilon}  1_{   [0 , \varepsilon]   }   (x)
  =
  \hat p (x) \int_W p(t,x) \, 1_{   [0 , \varepsilon]   }    (\hat \tau_{t} (x)) \nu(dt)
  + \left( 1 - \hat p (x) \right) 1_{   [0 ,   \varepsilon]   } (x)
   \quad \mbox{ for any $x \in [0,1]$} ,
$$
it follows from the definition of $ \hat \tau_t $ that 
 \begin{eqnarray*}
&  &  \hspace{-5ex}
 1_{   [0 ,  c)   } (\hat \tau_{t_0} (x)) \,  \left(  {\cal U}_{  \Upsilon} 
   1_{   [0 , \varepsilon]   }   \right) (\hat \tau_{t_0}(x)) \\
&=&
  \hat p (T_{t_0} (x) ) \int_W p(t, T_{t_0} (x)) \, 1_{   [0 , \varepsilon]   }    (\tau_2 (T_{t_0}(x))) \nu(dt)
  + \left( 1 - \hat p (T_{t_0} (x)) \right) 1_{   [0 ,   \varepsilon]   } (T_{t_0} (x)) \\
&=& 
  \hat p (T_{t_0} (x) )  1_{   [0 , \tau_2^{-1} (\varepsilon) ]   }   (T_{t_0}(x))
   + \left( 1 - \hat p (T_{t_0} (x)) \right) 1_{   [0 ,   \varepsilon]   } (T_{t_0} (x))
   \qquad \mbox{for any $x \in [c,1] $} .
\end{eqnarray*}

Thus,  in  (i),
setting $f =   {\cal U}_{  \Upsilon}    1_{   [0 , \varepsilon] }  $, 
we obtain (ii).
  \hfill $\Box$

\vspace{4ex}

\noindent 
{\it Proof of Theorem 6.5.}

Let $A$, $ \mu_A$, ${\cal U}_T$ and  $ \widetilde{{\cal U}}_T $ as in the proof of Theorem 6.2.
Since $h \le   \gamma_2$ on $ \bigcup_{t \in W} T_t|_{[c,1]}^{-1}(B(0)) $, it follows from Theorem 3.3 that  
\begin{eqnarray}
\mu ( [0 , \varepsilon] )
=   \sum_{n = 0}^{\infty}    \int_A
             { \widetilde{{\cal U}}}_T^n  \left(  {\cal U}_T 1_{ [ 0 , \varepsilon] }   \right) (x) \  \mu_A (dx) 
 \le    \gamma_2  \sum_{n = 0}^{\infty}    \int_A
             { \widetilde{{\cal U}}}_T^n  \left(  {\cal U}_T 1_{ [0  , \varepsilon] }   \right) (x) \  {m}|_A (dx).
\end{eqnarray}
For any small $    \varepsilon  \in (0, c_{\ast}) $
it follows from the assumption of the theorem that 
\begin{eqnarray}
&  &  \hspace{-5ex} \int_A
             { \widetilde{{\cal U}}}_T  \left(  {\cal U}_T 1_{   [0 , \varepsilon]   }   \right) (x) \  {m}|_A (dx)  
             \nonumber \\
&  = &
 \int_{W^2}  \left( 
         \int_{A}   
        p(t_0, x)  \, 1_{[ 0 ,  {T_{t_1} |_{[0, c) }}^{-1} (\varepsilon) ] } 
                         ( T_{t_0} (x)) \  p(t_1, T_{t_0} (x))   \, m(dx)  
         \right) \nu (dt_0) \nu (dt_1)    \nonumber 
\\ 
& = &  \int_{A}  \int_W \left( 
         \int_{W_2(\tau_2)} p(t_0, x)  \, 1_{[ 0 ,  {T_{t_1} |_{[0, c) }}^{-1} (\varepsilon) ] } 
                         ( T_{t_0} (x)) \  p(t_1, T_{t_0} (x))   \nu (dt_1) \right.  \nonumber  \\
& &  \qquad    +   \left.  
     \int_{W \setminus W_2(\tau_2)} p(t_0, x)  \, 1_{[ 0 ,  {T_{t_1} |_{[0, c) }}^{-1} (\varepsilon) ] } 
                         ( T_{t_0} (x)) \  p(t_1, T_{t_0} (x))  \nu (dt_1) 
     \right) \nu(dt_0) m(dx)      \nonumber  \\      
 &\le &  \int_{A }  \int_W \left( 
         \int_{W_2(\tau_2)} p(t_0, x)  \, 1_{ [ 0  ,  \tau_2 ^{-1}(\varepsilon) ]  } 
                         ( T_{t_0} (x)) \  p(t_1, T_{t_0} (x))   \nu (dt_1) \right.  \nonumber  \\
& &  \qquad    +   \left.  
     \int_{W \setminus W_2(\tau_2)} p(t_0, x)  \, 1_{[ 0,   \varepsilon]} 
                         ( T_{t_0} (x)) \  p(t_1, T_{t_0} (x))  \nu (dt_1) 
     \right) \nu(dt_0) m(dx)          \nonumber      \\
 & =  & 
  \int_{A }  \int_W p(t_0, x)       \nonumber      \\
  & &  \qquad 
  \left(  1_{ [0  ,  \tau_2^{-1}(\varepsilon) ]  }   ( T_{t_0} (x)) 
       \hat{p}  ( T_{t_0} (x))    + 
      1_{[ 0 ,   \varepsilon]} 
                         ( T_{t_0} (x))   (1 -   \hat{p}( T_{t_0} (x)) )       \right) \nu(dt_0) m(dx) .    
\end{eqnarray}
By Lemma 6.9(ii),  
 the right hand of the inequality (6.5) is 
$$
 \int_A
 { \widetilde{{\cal U}}}_{  \Upsilon}   \left(  {\cal U}_{  \Upsilon} 
   1_{   [0 , \varepsilon]   }   \right) (x) \  {m}|_A (dx). 
$$
Thus, we obtain 
$$
 \int_A { \widetilde{{\cal U}}}_T  \left(  {\cal U}_T 1_{   [0 , \varepsilon]   }   \right) (x) \  {m}|_A (dx)  
\le
 \int_A
 { \widetilde{{\cal U}}}_{  \Upsilon}   \left(  {\cal U}_{  \Upsilon} 
   1_{   [0 , \varepsilon]   }   \right) (x) \  {m}|_A (dx). 
$$
Similarly we can obtain
$$
 \int_A { \widetilde{{\cal U}}}_T^n  \left(  {\cal U}_T  1_{   [0 , \varepsilon]   }   \right) (x) \  {m}|_A (dx)  
\le
 \int_A
 { \widetilde{{\cal U}}}_{  \Upsilon}^n   \left(  {\cal U}_{  \Upsilon} 
   1_{   [0 , \varepsilon]   }   \right) (x) \  {m}|_A (dx)
$$
and 
\begin{equation}
 \sum_{n=0}^{\infty} \int_A { \widetilde{{\cal U}}}_T^n  \left(  {\cal U}_T 1_{   [0 , \varepsilon]   }   \right) (x) \  {m}|_A (dx)  
\le
 \sum_{n=0}^{\infty}
 \int_A
 { \widetilde{{\cal U}}}_{  \Upsilon}^n   \left(  {\cal U}_{  \Upsilon} 
   1_{   [0 , \varepsilon]   }   \right) (x) \  {m}|_A (dx).
\end{equation}
Since $\hat h = \hat {\hat h} $ on $A$, 
we have
\begin{eqnarray}
 \sum_{n = 0}^{\infty}    \int_A
             { \widetilde{{\cal U}}}_{  \Upsilon}^n   \left(  {\cal U}_{  \Upsilon} 
   1_{   [0 , \varepsilon]   }   \right) (x) \  {m}|_A (dx)
 & \le &   \frac 1 {\hat \gamma_1 }   
  \sum_{n = 0}^{\infty}    \int_A
           { \widetilde{{\cal U}}}_{  \Upsilon}^n   \left(  {\cal U}_{  \Upsilon} 
   1_{   [0 , \varepsilon]   }   \right) (x) 
             \,   \hat{\mu}|_A (dx)  
 \nonumber  \\           
  & = &  \frac 1 {\hat \gamma_1 }    \hat{\hat \mu} ( [0 , \varepsilon] )
\end{eqnarray}

Therefore, by (6.4), (6.6), (6.7) and Lemma 6.8, we obtain
$$
\mu([0, \varepsilon ]) 
\leq  
 \frac {\gamma_2 }{\hat  \gamma_1  \inf_{ x \in [0, c_\ast)} \hat p (x)} 
\hat{\mu}( [0, \varepsilon]  ) \quad
\mbox{for any $  \varepsilon  \in (0, c_{\ast})$.} 
$$
which shows the theorem. 
  \hfill $\Box$

\vspace{2ex}

\noindent 
{\bf 
Remark 6.10.} Under  the  assumption of Theorems 6.2 and  6.5, 
if $p(t,x) $ in the random map $T$ does not depend on $x$,
by a minor modification of the proofs of the theorems 
 we can show that 
$$
  \frac{\gamma_1 }{ \bar \gamma_2 } \bar{\mu}( [ \varepsilon_1 , \varepsilon_2 ]  ) 
  \leq  
 \mu([ \varepsilon_1 , \varepsilon_2 ]) 
 \leq  
 \frac {\gamma_2 }{\hat  \gamma_1  \inf_{ x \in [0, c_\ast)} \hat p (x)} 
\hat{\mu}( [\varepsilon_1 , \varepsilon_2 ]  ) \quad
\mbox{for any small $   0 < \varepsilon_1 <  \varepsilon_2 $.}
$$

\vspace{3ex}

\noindent 
{\bf 
Remark 6.11. }
If $T_t $'s  on $[c, 1] $ are  common, then 
$\bar \tau$ and $\hat \tau $ are deterministic maps.
Thus, in such a case, 
we can relatively easily apply our comparison theorems and 
we can estimate $ \mu([0, \varepsilon])$ 
by using invariant measures of the two deterministic maps

\vspace{2ex}

\section{How to apply the comparison theorems}

We will show how to apply our comparison theorems to some random dynamical systems. 
Here,  we mainly consider random maps with
 unbounded derivatives and/or 
 with an indifferent fixed point. 
 The class of random maps in this section includes Examples 1.1 and 1.2 in Section 1 
 and also includes  the following example.

\vspace{2ex}

\noindent
{\bf Example 7.1. }
Let $a_0$ and $b_0$ be constants with $ 1 < a_0 < b_0$ and 
let $a_1$ and $b_1$ be constants with $ 1 < a_1 < b_1 < 3$. 
Define $ T_{t,s}: [0,1] \to [0,1]$ by 
$$
 T_{t,s} (x) =  
 \left\{
\begin{array}{ll} 
x + 4^{-1} (4/3)^t x^t  \  &   \mbox{ for }  x  \in  [0 , 3/4),     
\\ 
\frac s {s+1} \left(  4x -3  \right)^{1/s}    &    \mbox{ for }  x  \in  [3/4, 1 ],
\end{array}  \right.
$$
where $(t, s)$ is randomly selected from $W =  [a_0, b_0] \times [a_1, b_1] $ accordingly to
the probability density function $ p(t, s, x) = 1/ \{(b_0 - a_0)( b_1 - a_1) \}$.

\vspace{2ex}

In Example 7.1, the partition of $[0, 1] $  is fixed, though the partitions of $[0, 1] $ in Examples 1.1 and 1.2 vary with $t$.
On the other hand,
$ T_{t,s} ( [3/4, 1 ] ) \ne [0, 1]$ in Example 7.1,  
though $ T_{t,s}$   in Examples 1.1 and 1.2 are  piecewise onto maps except an endpoint.
By some technical reasons,  we assume that each map $ T_{t,s}$ is  piecewise onto 
if the partition of $ [0,1] $ varies with the parameters.

\vspace{1ex}

Generalizing  Examples 1.1, 1.2 and 7.1, 
we consider the following setting.

\vspace{1ex}

Let $W \subset {\mathbb R} ^2$ with a  measure $\nu = \nu_0 \times \nu_1$. 
(In Examples 1.1, 1.2 and 7.1,  $\nu$ is the Lebesgue measure.)
Put 
$
V_0(s) = \{ t : (t,s) \in W \}
$
and 
$
V_1(t) = \{ s : (t,s) \in W \}
$.
For a probability density function $ p((t, s), x)$
such that $ \int_W p((t, s), x) \nu_0(dt ) \nu_1 (ds)  = 1 $ for each $ x \in [0,1]$, 
put 
$$
 p_0(t,x) = \int_{ V_1(t) }  p((t, s), x) \nu_1 (ds)
 \quad \mbox{and} \quad
 p_1(s,x) = \int_{ V_0 (s) }  p((t, s), x) \nu_0 (dt). 
$$
We consider  piecewise $ C^1 $ random maps
$ T = \{ T_{t, s} ; p((t, s) , x); \{[0, c_t), [c_t, 1] \} : {(t, s)} \in W \} $ 
 which satisfy the following assumption (7A).

\noindent 

\begin{itemize}
\item[{(7A) } ]

\begin{itemize}
\item [(i)]
$ a_0 = \inf \{t : (t, s) \in W \}$ and 
$b_0 = \sup  \{t : (t, s) \in W \}$ 
satisfy $ 1 < a_0  \le b_0 < \infty$,  
and 
 $a_1 = \inf \{s : (t, s) \in W \}$ and 
 $ b_1= \sup  \{s  : (t, s) \in W \}$ satisfy $1 \le  a_1 \le  b_1 < \infty $.  

\item [(ii)]
For each  $ (t, s) \in W $,  $T_{t, s}(0) = 0$ . Put $c = \inf_{  (t, s) \in W } c_t > 0$. Further, 

(a)
If $c_t =c$ for all $t$, we assume $T_{t,s} ([0, c) ) = [0, T_{t,s}(c) ) $ 
and  $T_{t,s} ([c, 1] ) = [0, T_{t,s}(1) ]$.

(b) 
If $c_t $ depends on $t$, we assume $T_{t,s} ([0, c_t) ) = [0, 1 )$ 
and  $T_{t,s} ([c_t , 1] ) = [0, 1 ]$.

\item [(iii)]
There exists a constant $M_0$ such that 
$\displaystyle 
 \bigvee_{[0,1]} \frac { p((t, s), \cdot )}{T_{t,s}^\prime}
 < M_0 $ 
for a.s. $(t,s) \in W$.

\item [(iv)] If $T_{t,s}^\prime (0) = 1$ for  a.s. $(t,s) \in W$,   
$\displaystyle \frac { p((t, s), x )}{T_{t,s}^\prime (x) } $ is monotonically  increasing as $x$ goes to $0$ for each $ (t, s) \in W$.

\item [(v)]
$
T_{t,s}(x) =\ell (t, s) (x-c_t)^{1/s}  
\  \mbox{on } T_{t,s}|_{ [c_t , 1]}^{-1} B(0), 
$
where 
 $\ell (t, s)  $ is a function of $(t, s) \in W  $ such  that 
$ M_1  \leq (\ell ( t, s) )^s \leq M_2$ with  positive constants 
$M_1$ and $M_2$.

\item [(vi)]
$
T_{t,s}^{\prime}(x) > 1 
$
on $ (c_t, 1] $ for each $ (t,s) \in W$.

\end{itemize}

\end{itemize}

\vspace{2ex}

Using our comparison theorems, 
for random maps with the indifferent fixed point 0 and unbounded derivatives, 
we will show the following two theorems.

\vspace{3ex}

\noindent 
{\bf Theorem 7.1. }
{\em 
Let 
$ T = \{ T_{\bm t}  ; p(\bm t, x) : \bm t = (t, s)  \in W \} $ be a piecewise $ C^1 $ random map
which satisfies the assumption (7A). 
Assume that 
for  a.s. $\bm t \in W$,  
 $$  T_{\bm t}(x)  =    x  +  m_t \, x^{t} +  o( x^{t} )  
 \   \mbox{  in } \  [0,  \varepsilon_0 ),  
 $$
where $  \varepsilon_0  > 0 $ is a constant which is independent of $\bm t$ 
and $m_t >0 $ is a constant for each $t $ with $(t,s) \in W $. 
Let $ \nu_1$ be 
the Lebesgue measure. 
Assume that  $p( \bm t, x) $ satisfies 
 the following condition: 
 
(a) there exists a set of constants
$x_0 , x_1 \in [c,1]$,  $C_1 >0 $ and $C_2 >0$ such that 
$$
C_1 p(\bm t, x_0)  \le p(\bm t, x) \le C_2 p(\bm t, x_1)  \quad \mbox{for any} \quad  x \in [c, 1].
$$

(b) $p_1(s, x_0) \ge \bar p_1$ and 
$p_1(s, x_1) \le  p_1^{\ast}$, 
where $\bar p_1 >0$ and $p_1^{\ast} >0 $ are constants.

(c) For some   small  constant $ \kappa \geq 0 $, 
$ \displaystyle  \inf_{x \in [0,  c]} \int_{[a_0, a_0 + \kappa ] }  p_0(t, x) \nu_0 (dt) > 0$.

Then, 
the random map $T$ has 
an  absolutely continuous $ \sigma$-finite invariant measure $ \mu$
which satisfies the following (1) and (2) : 
\begin{itemize}
\item[(1)]
If  $ a_0 \ge a_1 +1$, then 
$ \mu([0, \varepsilon]) = \infty $ for any small $ \varepsilon > 0$.

\item[(2)]
If  $ a_0  <   a_1 +1 $, then
$$ 
\frac {C_3  \varepsilon^{  a_1 +1 - a_0  } }  {- \log \varepsilon } \le  \mu([0, \varepsilon]) 
 \le  \frac { C_4  \varepsilon^{  a_1  +1 - a_0 - \kappa } }  {- \log \varepsilon }
\quad \mbox{ for any small $ \varepsilon > 0$,}
$$
where $C_3$ and $C_4$ are constants.

\end{itemize}

}

\vspace{3ex}

\noindent  
 {\bf Theorem 7.2.}
 {\em 
 \ Assume that $\cup_t  V_1(t) $ is a finite set. 
In Theorem 7.1,  if we replace `$ \nu_1$ is the Lebesgue measure' by 
`$\nu_1$ is  the counting measure', 
 then the conclusion (1) of Theorem 7.1 remains valid.
 However, under this situation, 
 the conclusion (2) of Theorem 7.1 is changed to the following:
\begin{itemize}
\item[(2)$ ^\ast$]  If  $ a_0  <   a_1 +1 $, then
$$ 
C_3^{\ast}  \varepsilon^{  a_1 +1 - a_0  } \le  \mu([0, \varepsilon]) 
 \le  C_4^{\ast}  \varepsilon^{  a_1  +1 - a_0 - \kappa } 
\quad \mbox{ for any small $ \varepsilon > 0$,}
$$
where $C_3^{\ast}$ and $C_4^{\ast}$ are constants.

\end{itemize} 
 }

\vspace{2ex}

\noindent 
{\bf Remark 7.3.} 
In Theorems 7.1 and 7.2,
if $a_0$ is selected with a positive probability, 
the condition (c) is satisfied for $ \kappa = 0$.
Thus, $ \kappa = 0$ 
in the estimates of  (2) and  (2)$ ^\ast$
under such a  situation.

\vspace{2ex}

Theorem 7.2 can be adapted to deterministic maps. By Theorem 7.2 we obtain:

\vspace{2ex}

\noindent 
{\bf Corollary  7.4.}  {\em 
Let $T$ be  a deterministic map which satisfies the assumption  of Theorem 7.2 with $W = \{ (a_0, a_1) \}$. 
Then $T$ has an absolutely continuous $ \sigma$-finite invariant measure $\mu$ such that 
\begin{eqnarray*}
& &  \mu([0, \varepsilon]) = \infty \quad \mbox{for any small $ \varepsilon > 0$ if $ a_0 \ge a_1 +1$, } \\
& & 
 \mu([0, \varepsilon]) 
 \approx   \varepsilon^{  a_1  +1 - a_0 } 
\mbox{ if $ a_0  < a_1 +1$. }
\end{eqnarray*}

}

 \vspace{2ex}

From Corollary  7.4, each individual map $ T_{t,s}$ in Example 1.1
 has an absolutely continuous $ \sigma$-finite invariant measure $\mu$ 
such that 
\begin{eqnarray*}
& &  \mu([0, \varepsilon]) = \infty \quad \mbox{for any small $ \varepsilon > 0$ if $ t \ge s +1$, } \\
& &   \mu([0, \varepsilon]) 
\approx  \varepsilon^{  s +1 - t } 
\quad \mbox{  if  $t  <  s +1$}. 
\end{eqnarray*} 
 The random map $T$ of  Example 1.1  satisfies the assumption of  Theorem 7.1,  
 by setting  $a_1 = a_0 c_0 $. 
 The random map $T$ of  Example 7.1  also satisfies the assumption of  Theorem 7.1.  
In these cases,  the constant  $ \kappa > 0$ can be chosen arbitrary small.

 \vspace{1ex}

 Applying our comparison theorems 
to random maps with  unbounded derivatives and without indifferent fixed points,  
we will also show the following two theorems.

\vspace{3ex}

\noindent 
{\bf Theorem 7.5. }
{ \em 
Let 
$ T = \{ T_{\bm t } ; p( \bm t, x) : \bm t \in W \} $ be a piecewise $ C^1 $ random map
which satisfies the  assumption (7A).  
Assume that 
for  a.s. $\bm t \in W$,  
 $$  x < a_0 x \le T_{\bm t}(x)  \le b_0 x
 \   \mbox{  in } \  [0,  \varepsilon_0 ),  
 $$
where $  \varepsilon_0  > 0 $ is a constant which is independent of $\bm t$. 
Let  $\nu_1$ be the Lebesgue measure and let $ p( \bm t , x)$ be as in Theorem 7.1.

Then, 
the random map $T$ has 
an  absolutely continuous invariant probability measure $ \mu$ 
which satisfies 
$
\mu([0, \varepsilon])
\approx    \varepsilon^{a_1 } / (- \log \varepsilon )  .
$
}

\vspace{4ex}

\noindent 
{\bf Theorem 7.6. }
{\em 
  Assume that $\cup_t  V_1(t) $ is a finite set. 
In Theorem 7.5, if we replace `$ \nu_1$ is the Lebesgue measure' by 
`$\nu_1$ is  the counting measure', 
 then 
 the random map $T$ has 
an  absolutely continuous invariant probability measure $ \mu$ 
which satisfies 
$
  \mu([0, \varepsilon]) \approx  \varepsilon^{ a_1} .
$

}

\vspace{3ex}

Theorem 7.6  can be adapted to deterministic maps.  

\vspace{2ex}

\noindent 
{\bf Corollary  7.7. }
{\em 
Let $T$ be a deterministic map which satisfies the assumption of Theorem 7.6 with 
$\cup_t  V_1(t)  = \{  a_1\}$. 
Then,  the absolutely continuous invariant probability measure $ \mu$ satisfies  
$
  \mu([0, \varepsilon]) \approx  \varepsilon^{ a_1} .
$
}

\vspace{3ex}

From Corollary  7.7, each individual map $ T_{t,s}$ in Example 1.2  
has an absolutely continuous invariant probability measure $\mu$ 
which satisfies  
$
  \mu([0, \varepsilon]) \approx  \varepsilon^{ s } .
$
 The random map $T$ of  Example 1.2  satisfies the assumption of  Theorem 7.5,  
 setting  $a_1 = a_0 c_0 $.

 \vspace{3ex}

Since, in our comparison theorems, we assume that $T$ satisfies (6A),  it is important that 
the first return random map of $T$ 
has  an absolutely continuous invariant measure whose density 
is  bounded above and bounded away from 0 on $  \bigcup_{ \bm t \in W} T_{\bm t} |_{[c,1]}^{-1}(B(0)) $. 
Since the random map $T$ as in Theorem 7.1  satisfies the assumption of Theorem 5.1, 
 the random map $T$ has the absolutely continuous invariant measure whose density is 
bounded above on $ [c, 1]$.  
To study the lower estimates of the invariant density of $T$ and of its auxiliary random maps, 
we prepare some lemmas.

Let $\tau_0 : [0, c ) \to [0, 1) $ be a $C^1$ and surjective map 
such that 
 $ \tau_0^{\prime} (x) > 1$ on $ (0, c)$.
For each $\bm t \in W$,
define $ \tilde \tau_{ \bm t} : [0, 1]  \to [0,1]$ 
by 

$$
{ \tilde \tau}_{\bm t} (x) := 
\left\{
\begin{array}{ll}  
\tau_0 (x),   &    \quad  [0, c),  \\   
T_{\bm t}(x), &  \quad  [ c, 1].
\end{array}
 \right.
$$
In the following three lemmas, 
we also consider the random map 
$ \tilde \tau =  \{ \tilde \tau_{\bm t} ; p(\bm t, x) : \bm t \in W \} $
to study the auxiliary random maps.  

 \vspace{4ex}

\noindent 
{\bf Lemma 7.8. } 
{\em 
  Let $T$ be a random map as in (7A) with (a) of (ii). 
 Assume that  $ c  = c_{t}$ is independent of $ (t, s) \in W$.
Let $ T_{\bm t}^{\ast}$ be a continuous extension  of $T_{\bm t} |_{[0, c)}$ to $ [0, c]$ 
and put  
$ c^\ast = \inf_{ \bm t \in W}  (T_{\bm t}^{\ast} (  c  )  ) $. 
Consider the first return random map   
$ R= \{ r_u ; Pr(u, x) : u \in  \widetilde{W}\} $
of $T$ on $[ c  ,  c^\ast]  $.
Let $ J  \subset [c ,  c^\ast]$ be a nontrivial  interval.  
Then,  there is an $n$
such that  
$$
\tilde{\nu}^n \left\{ \displaystyle (u_1, u_2,  \ldots,  u_n) \in  {\widetilde{W}}^n  :    
r_{u_n}  \circ  \cdots  \circ r_{u_2}   \circ r_{u_1}   ( J)   \supset [c,  c^\ast ]    \right\} > 0.
$$
Further, for the random map $\tilde \tau $,  we have the same conclusion.  
}

\vspace{3ex}

\noindent 
{\em Proof.}
Let ${\cal P}_u  $ be the partition corresponding to $ r_u $.
If a nontrivial interval $J$ is contained in one of the intervals in  ${\cal P}_u  $,
then $ r_{u} (J )$ is an interval and $ m(r_u(J)) > m(J)  $.
Thus,
 there is an $N$  such that $ r_{u_N}  \circ  \cdots  \circ r_{u_2}   \circ r_{u_1}   ( J) $
contains an endpoint of one of the intervals in ${\cal P}_{u_{N +1 }}  $.
This implies that 
$$
\tilde{\nu}^{N+1} \left\{ \displaystyle (u_1, u_2,  \ldots,  u_{N+1}) \in  {\widetilde{W}}^{N+1}  :    
r_{u_{N+1}}  \circ  \cdots  \circ r_{u_2}   \circ r_{u_1}   ( J)  
      \supset [c,  c +   \epsilon_0 ]    \right\} > 0
$$
for a small $  \epsilon_0 > 0 $.
Since  $r_u ( [c,  c +   \epsilon_0 ]  )  \supset  [c, c^\ast]$, 
we have 
$$
\tilde{\nu}^{N+2} \left\{ \displaystyle (u_1, u_2,  \ldots,  u_{N+2}) \in  {\widetilde{W}}^{N+2}  :    
r_{u_{N+2}}  \circ  \cdots  \circ r_{u_2}   \circ r_{u_1}   ( J)  
      \supset [c,  c^\ast ]    \right\} > 0.
$$
Therefore, we obtain the lemma for $T$.   
Similarly,  we obtain the lemma for $\tilde \tau $.   
\hfill $\Box$

\vspace{3ex} 

\noindent 
{\bf Lemma 7.9.}  
{\em 
 Let $T$ be a random map as in  (7A) with (b) of (ii). 
Consider the first return random map   
$ R= \{ r_u ; Pr(u, x) : u \in  \widetilde{W}\} $
of $T$ on $[c  ,  1]  $.
Let $ J  \subset [c ,  1]$ be a nontrivial  interval.  
Then,  there is an $n$  such that 
$$ 
\tilde{\nu}^n \left\{ \displaystyle (u_1, u_2,  \ldots,  u_n) \in  {\widetilde{W}}^n  :    
r_{u_n}  \circ  \cdots  \circ r_{u_2}   \circ r_{u_1}   ( J)   \supset [c, 1 ]   \right\} > 0.
$$
Further, for the random map $\tilde \tau $,  we have the same conclusion.

}

\vspace{1ex}

\noindent 
{\em Proof.} 
 In the same way to the proof of the previous lemma,  
 there is an $N_0$  such that $ r_{u_{N_0}}  \circ  \cdots  \circ r_{u_2}   \circ r_{u_1}   ( J) $
contains an endpoint of one of the intervals in ${\cal P}_{u_{{N_0 +1}}}  $.
This implies that there exists an $ \epsilon_0 > 0 $ such that 
$$
\tilde{\nu}^{N_0+1} \left\{ \displaystyle (u_1, u_2,  \ldots,  u_{N_0+1}) \in  {\widetilde{W}}^{N_0+1}  :    
r_{u_{{N_0 +1} } }  \circ  \cdots  \circ r_{u_2}   \circ r_{u_1}   ( J)  \supset [1 - \epsilon_0 ,  1]    \right\} > 0. 
$$
Since there exists an $ N_1 $ such that 
$ r_{u_{N_1}}  \circ  \cdots  \circ r_{u_{N_0 +2}  }      ( [1 - \epsilon_0 ,  1]) 
    \supset [c,  1] $,
we obtain 
$$
\tilde{\nu}^{N_1} \left\{ \displaystyle (u_1, u_2,  \ldots,  u_{N_1}) \in  {\widetilde{W}}^{N_1}  :    
r_{u_{{N_1} } }  \circ  \cdots  \circ r_{u_2}   \circ r_{u_1}   ( J)  \supset [c ,  1]    \right\} > 0. 
$$
Therefore, we obtain the lemma for $T$.   
Similarly,  we obtain the lemma for $\tilde \tau $.  \hfill $\Box$

\vspace{3ex}

\noindent 
{\bf  Lemma 7.10.} 
{\em   
 Let $T$ be a random map as in (7A). 
Let $C \subset A $ be a closed interval or a union of closed intervals.  
Let $ R= \{ r_u; Pr(u, x) : u \in  \widetilde{W} \} $ 
be the first return random map of $T$ on $C $   and 
let  ${\cal L }_R  $ be the Perron-Frobenius operator corresponding to $R$.
Let $ J  \subset C$ be a nontrivial  interval. 
Assume that there is an $n$
such that  
$$
\tilde{\nu}^n \left\{ \displaystyle (u_1, u_2,  \ldots,  u_n) \in  {\widetilde{W}}^n  :    
r_{u_n}  \circ  \cdots  \circ r_{u_2}   \circ r_{u_1}   ( J)   \supset C   \right\} > 0.
$$
Then,  there exists  a  function $h \ge 0$ of  bounded variation   
such that 
 $$ \int_{C} h(x) m(dx) =1,  \  \    
{\cal L }_R h = h ,   \  \   \mbox{ and}   \   \      \inf_{x \in C} h(x) > 0. $$ 
Further, for the random map $\tilde \tau $,  we have the same conclusion.

}
 
 \vspace{2ex}

\noindent 
{\em Proof.} 
By Theorem 5.1,  $R$ has a stationary density $h$ of bounded variation.
Hence, 
there exists a nontrivial interval $J \subset C $ and a constant $ \gamma > 0  $ 
such that $ h(x) \ge \gamma 1_J (x) $.  
Then, 
$$
h(x)  =  {\cal L }_R^n   h(x)  \ge \gamma \  {\cal L }_R^n  1_J (x) 
\quad  \mbox{for} \  x \in C.
$$
Put 
$$
{\widetilde{W}}^n_{\ast} =  \left\{ \displaystyle (u_1, u_2,  \ldots,  u_n) \in  {\widetilde{W}}^n  :    
r_{u_n}  \circ  \cdots  \circ r_{u_2}   \circ r_{u_1}   ( J)   \supset C   \right\}.
$$
By the assumption, there exists  an $n$ with $  \tilde{\nu}^n( \widetilde{W}_{\ast}^n  ) > 0$.   
For this $n$
\begin{eqnarray*}
{\cal L }_R^{ n}  1_J (x) 
&=&
 \int_{ {\widetilde{W}}^n}  \sum_{y \in (r_{u_n} \circ  \cdots  \circ r_{u_1} )^{-1} \{ x \}  }
\frac{Pr(u_1,  \ldots  ,  u_n ,y)}
       {(r_{u_n} \circ  \cdots  \circ r_{u_1} )^{\prime} (y)}
       1_J (y)   \tilde{\nu}(du_1) \cdots \tilde{\nu}(du_n)  \\
&   \ge &
 \int_{{\widetilde{W}}^n_{\ast}}  \sum_{y \in (r_{u_n} \circ  \cdots  \circ r_{u_1} )^{-1} \{ x \}  }
\frac{Pr( u_1,  \ldots , u_n ,y)}
       {(r_{u_n} \circ  \cdots  \circ r_{u_1} )^{\prime} (y)}
       1_J (y)   \tilde{\nu}(du_1) \cdots  \tilde{\nu}(du_n),        
\end{eqnarray*}
where 
$  Pr(u_1,  \ldots , u_n ,y) = Pr(u_1, y ) \cdots  Pr (u_n,  r_{u_{n-1}} \circ \cdots \circ r_{ u_{1}} (y) )$.
By the definition of $ {\widetilde{W}}^n_{\ast} $, 
if  $ (u_1, u_2,  \ldots,  u_n) \in {\widetilde{W}}^n_{\ast}$, 
$  (r_{u_n} \circ  \cdots  \circ r_{u_1} )^{-1} \{ x \}  \cap J \ne \emptyset $ for 
$ x \in  C $.

Thus,   $  {\cal L }_R^n  1_J (x)  > 0$  on $ C $.
Therefore, $ h(x) > 0 $ on $  C $, 
which shows the lemma for $T$.   
Similarly,  we obtain the lemma for $\tilde \tau $.
  \hfill $\Box$

\vspace{5ex}

As in the previous section, 
in the proofs of the Theorems 7.1, 7.2, 7.5 and 7.6, put $B(0) = [0, \varepsilon_0] $. 

\vspace{3ex}

\noindent 
{\it Proof of Theorem 7.1.}
First we consider the lower estimate and after that we consider the upper estimate.

\noindent 
[Lower estimate]

Let $\bar  c$ be a small constant with $0 < \bar c < \min \{ c,   \varepsilon_0 \}$.
Dfine $\bar m = \sup \{ m_t : \   (t, s) \in W  \} $. \ 
Let $\tau_1 :[0, c ) \to [0, 1) $ be a $C^1$ and surjective map 
such that 
$\tau_1(x) = x +  \bar m x^{a_0} $ on $ [0,   \bar c ] $ and that $ \tau_1^{\prime} (x) > 1$ on $ (0, c)$.
Then, for a.s. $ \bm t \in W$, $ T_{ \bm t }  (0,  \varepsilon ) \subset \tau_1 (0,  \varepsilon )$ 
for any $ \varepsilon \in (0 , \bar c )$. 
For each $\bm t \in W$,
define $ \bar{\tau}_{\bm t} : [0, 1]  \to [0,1]$ 
by 
$$
\bar{\tau}_{\bm t} (x) := 
\left\{
\begin{array}{ll}  
\tau_1 (x),   &    \quad  [0, c),  \\   
T_{\bm t}(x), &  \quad  [ c, 1].
\end{array}
 \right.
$$
Then, by Theorem 5.1 
the random map 
$ \bar{\tau} = \{ \bar {\tau}_{\bm t} , p(\bm t,x) : \bm t \in W \}$
has an absolutely continuous invariant measure $ \bar {\mu}$ with $\bar \mu([c,1]) =1 $.
Let $\bar h$ be a density of $ \bar {\mu}$.
Then, 
by Theorem 5.1 and Lemmas 7.8, 7.9 and 7.10,  we have constants $\bar \gamma_1 $ and $\bar \gamma_2$ such that 
$ 0<  \bar \gamma_1 \le  \bar{h} (x) \le  \bar \gamma_2 < \infty$ on $  \bigcup_{ \bm t \in W} T_{\bm t} |_{[c,1]}^{-1}(B(0)) $.   

Put $ A= [c, 1]$ and let $ \bar \mu_A = \bar \mu |_A $
be the invariant probability measure of the first return map of the random map $ \bar{\tau} $.
Then, 
by Theorem 3.3,  we have 
\begin{equation}
\bar \mu ([0, \varepsilon]) = 
\sum_{n=0}^{\infty} 
  \int_{ A }
             { \widetilde{{\cal U}}^n_{\bar \tau} }  \left(  {\cal U}_{\bar \tau}   1_{  [0, \varepsilon]  }   \right) (x) \  \bar\mu_A (dx). 
\end{equation}
Since $\bar h \ge \bar \gamma_1$ on $  \bigcup_{ \bm t \in W} T_{\bm t} |_{[c,1]}^{-1}(B(0)) $, 
we have 
\begin{eqnarray}
  \int_{ A }
             { \widetilde{{\cal U}}^n_{\bar \tau} }  \left(  {\cal U}_{\bar \tau}   1_{  [0, \varepsilon]  }   \right) (x) \  \bar \mu_A (dx) 
 &  = &             
    \int_{W}  \left( 
         \int_{A}   
        p({\bm t}, x)  \, 1_{ [ 0,  \tau_1^{-n} (\varepsilon) ] } 
                         ( T_{{\bm t}} (x))    \, \bar \mu_A (dx)  
         \right) \nu (d {\bm t})   \nonumber  \\
 &  \ge &   
\bar \gamma_1  
    \int_{W}  \left( 
         \int_{A}   
        p({\bm t}, x)  \, 1_{ [ 0,  \tau_1^{-n} (\varepsilon) ] } 
                         ( T_{{\bm t}} (x)) \,   m (dx)  
         \right) \nu (d {\bm t}).
\end{eqnarray}
for any  small $ \varepsilon >0 $. 
Since it follows from (7A)(v) that 
$$
0 \leq \ell (t, s) (x-c_t)^{1/s }  \leq \tau_1^{-n} (\varepsilon) 
\quad  \Longleftrightarrow  \quad 
c_t \leq x \leq c_t + \frac{ \left( \tau_1^{-n} (\varepsilon) \right)^s }{  \left( \ell (t, s)  \right)^s },
$$
by the condition (a) we obtain 
\begin{eqnarray*}
 \int_{W}   \int_{A}   
        p({\bm t}, x)  \, 1_{ [ 0,  \tau_1^{-n} (\varepsilon) ] } 
                         ( T_{{\bm t}} (x))  \, m(dx)   \nu (d {\bm t}) 
   &  \geq  &
  C_1   \int_{W}  
     p({\bm t}, x_0)     \frac{ \left( \tau_1^{-n} (\varepsilon) \right)^s }{M_2}  \nu (d {\bm t})  \\
 & = &
 C_1 \int_{[a_1, b_1]}  p_1(s, x_0 )  \frac{ \left( \tau_1^{-n} (\varepsilon) \right)^s }{M_2}  \nu_1(ds).       
\end{eqnarray*}
Put 
$
C_5 =  \bar p_1 C_1 / M_2. 
$
Since $ \nu_1$ is the Lebesgue measure on $[a_1, b_1] $, 
by the condition (b) we have 
\begin{eqnarray}
& &  \int_{W}   \int_{A}   
        p({\bm t}, x)  \, 1_{ [ 0,  \tau_1^{-n} (\varepsilon) ] } 
                         ( T_{{\bm t}} (x))  \, m(dx)   \nu (d {\bm t}) 
 \geq  C_5 \int_{a_1 }^{ b_1 }  \left( \tau_1^{-n} (\varepsilon) \right)^s ds
              \\
 &  = &   
\frac{C_5}{ \log  \tau_1^{-n} (\varepsilon)  } 
  \left( \left( \tau_1^{-n} (\varepsilon) \right)^{ b_1 } -  \left( \tau_1^{-n} (\varepsilon) \right)^{ a_1 } \right).   \nonumber 
\end{eqnarray}
It is easy to see that 
$
C_{\tau_{1,1}}  \left(  \varepsilon^{-(a_0 -1)} + n  \right)^{- 1/(a_0-1)} 
\le  \tau_1^{-n} (\varepsilon)
  \le C_{\tau_{1,2}}  \left(  \varepsilon^{-(a_0 -1)} + n  \right)^{- 1/(a_0-1)}  
$ for any small $\varepsilon >0 $, 
where $C_{\tau_{1,1}} > 0$ and $C_{\tau_{1,2}} > 1$ are constants.
Thus,  we have 
\begin{eqnarray}
 \int_{ a_1 }^{ b_1 }  \left( \tau_1^{-n} (\varepsilon) \right)^s ds
&\ge & 
\frac  { C_{\tau_{1,1}} \left(  \varepsilon^{-(a_0 -1)} + n  \right)^{- a_1 / (a_0 -1)} 
             -  C_{\tau_{1,2}} \left(  \varepsilon^{-(a_0 -1)} + n  \right)^{-{b_1} / (a_0 -1)} }
        {- \log C_{\tau_{1,2}}  + \frac 1 {a_0 - 1 } \log  \left(  \varepsilon^{-(a_0 -1)} + n  \right)} \nonumber \\
 &\ge &        
 \frac  { \frac 1 2 C_{\tau_{1,1}} \left(  \varepsilon^{-(a_0 -1)} + n  \right)^{- a_1  / (a_0 -1)} }
        { \frac 1 {a_0 - 1 } \log  \left(  \varepsilon^{-(a_0 -1)} + n  \right)}
\quad \mbox{for any small} \ \varepsilon >0.
\end{eqnarray}

(1) 
Assume $  a_0  \ge  a_1  +1$. 
Then,  we have $  a_1/ (a_0  -1)  \le 1 $.
Thus, by the inequalities (7.1)-(7.4) we have 
$
\bar \mu ([0, \varepsilon]) 
 = \infty 
$
for any $ \varepsilon > 0$. 

Therefore, by Theorem 6.2, we obtain 
$$
\mu ([0, \varepsilon]) = \infty  \quad \mbox {for any $ \varepsilon > 0$  \quad  if}  \quad    a_0 \ge  a_1 +1.
$$

(2)
Assume $  a_0  <  a_1  +1$. Then, by  the inequality (7.4) we have 
\begin{eqnarray}
\sum_{n=0}^{\infty} 
 \int_{ a_1 }^{ b_1 }  \left( \tau_1^{-n} (\varepsilon) \right)^s ds
& \ge & 
 \frac 1 2 \sum_{n=0}^{\infty} 
 \frac  {  C_{\tau_{1,1}} \left(  \varepsilon^{-(a_0 -1)} + n  \right)^{- a_1 / (a_0 -1)} }
        { \frac 1 {a_0 - 1 } \log  \left(  \varepsilon^{-(a_0 -1)} + n  \right)}  \nonumber \\
& \ge & 
 \frac 1 2 \sum_{n=0}^{ [  \varepsilon^{-(a_0 -1) }] } 
 \frac  { C_{\tau_{1,1}} \left(  \varepsilon^{-(a_0 -1)} + n  \right)^{- a_1 / (a_0 -1)} }
        { \frac 1 {a_0 - 1 } \log  \left(  \varepsilon^{-(a_0 -1)} + n  \right)} \nonumber  \\
& \ge & 
 \frac 1 2 \sum_{n=0}^{ [  \varepsilon^{-(a_0 -1) }]  } 
 \frac  {  C_{\tau_{1,1}} \left(  \varepsilon^{-(a_0 -1)} + n  \right)^{- a_1 / (a_0 -1)} }
        { \frac 1 {a_0 - 1 } \log  \left( 2  \varepsilon^{-(a_0 -1) }   \right)}  \nonumber \\        
& \ge &  \frac{  C_{\tau_{1,1}} (a_0 -1 ) ( 1 - 2^{ (a_0 -  a_1 -1)  / (a_0 - 1)} ) \varepsilon^{a_1 + 1 - a_0  } }
              {- 4 (a_1 + 1- a_0 ) \log  \varepsilon } 
\end{eqnarray}
for any   small $ \varepsilon >0$.
Put $ C_6 =  C_5  C_{\tau_{1,1}}  (a_0 -1 ) ( 1 - 2^{ (a_0 - a_1-1)  / (a_0 - 1)} ) / (4 (a_1+ 1 - a_0 ))$.
Then, by the inequalities (7.1), (7.2), (7.3) and (7.5)  we have
$$
\bar \mu ([0, \varepsilon]) 
 \ge \frac{  \bar \gamma_1 C_6 \  \varepsilon^{ a_1 +1 - a_0 } }{ - \log  \varepsilon  }
 \quad \mbox{for any  small $ \varepsilon >0$.}
$$
Therefore, by Theorem 6.2, we obtain 
$$
\mu ([0, \varepsilon]) 
\ge  \frac{ \gamma_1 \bar \gamma_1 C_6 \  \varepsilon^{a_1 + 1 - a_0  } }{ - \bar \gamma_2 \log  \varepsilon  }
 \quad \mbox{for any  small $ \varepsilon >0$}
  \quad \mbox {if}  \quad    a_0 <  a_1 +1, 
$$
which shows the lower estimate of (2).

\vspace{1ex}

\noindent 
[Upper estimate]

Let  $\hat  c$ be a small constant with  $0 < \hat c < \min \{ c,   \varepsilon_0 \}$.
Dfine $\hat m = \inf  \{ m_t : \  (t,s)  \in W \} $. \ 
Since $ a_0 < a_1 + 1$, 
we can choose  a  small constant $ \kappa > 0$  such that  $ a_0 + \kappa < a_1 + 1$. 
Let $\tau_2 :[0, c ) \to [0, 1) $ be a  $C^1$ and surjective map 
such that 
$\tau_2 (x) = x +  \hat m x^{a + \kappa} $ on  $ [0,   \hat c ] $ and that $ \tau_2^{\prime} (x) > 1$ on $ (0, c)$.
Put 
$$ W_2 (\tau_2) : = \{ \bm t \in W \ 
 | \  T_{ \bm t }  (0,  \varepsilon ) \supset \tau_2 (0,  \varepsilon ) \ \mbox{for any $ \varepsilon \in (0, \hat c )$}
\}
$$
 and 
$$
\hat p (x) :=
\int_{W_2(\tau_2)}  p( \bm t, x) \nu(d \bm t)
\quad \mbox{for $ x \in  [0, c) $}. 
$$
Then, by the condition (c) we have 
$
\hat p (x) > 0 
$
for $ x \in  [0, c) $.
For each $ \bm t \in W$,
define $ \hat{\tau}_{\bm t} : [0, 1]  \to [0,1]$ 
by 
$$
\hat{\tau}_{\bm t} (x) := 
\left\{
\begin{array}{ll}  
\tau_2 (x),   &    \quad  [0, c),  \\   
T_{\bm t} (x), &  \quad  [ c, 1].
\end{array}
 \right.
$$
By Theorem 5.1 the random map 
$ \hat{\tau} = \{ \hat{\tau}_{\bm t} , p(\bm t,x) ; \bm t \in W \}$
has an absolutely continuous invariant measure $ \hat{\mu}$
with $ \hat \mu([c,1]) =1 $.
Let $\hat h$ be a density of $ \hat{\mu}$.
As $ \bar h$ in the lower estimate, it follows  that  
$ 0< \hat \gamma_1 \le  \hat{h}(x)  \le  \hat \gamma_2 < \infty$
 on $  \bigcup_{\bm t \in W} T_{\bm t} |_{[c,1]}^{-1}(B(0)) $, 
 where $\hat \gamma_1 $ and $\hat \gamma_2 $ are constants.

Put $ A= [c, 1]$ and let $ \hat \mu_A = \hat \mu |_A $
be the invariant probability measure of the first return map of the random map $ \hat {\tau} $.
By Theorem 3.3,  we have 
\begin{equation}
\hat \mu ([0, \varepsilon]) = 
\sum_{n=0}^{\infty} 
  \int_{ A }
             { \widetilde{{\cal U}}^n_{\hat \tau} }  \left(  {\cal U}_{\hat \tau}   1_{  [0, \varepsilon]  }   \right) (x) \  \hat \mu_A (dx). 
\end{equation}
Since $  \hat{h}(x)  \le  \hat \gamma_2 $  on $  \bigcup_{\bm t \in W} T_{\bm t }|_{[c,1]}^{-1}(B(0)) $, 
we have 
\begin{equation}
  \int_{ A }
             { \widetilde{{\cal U}}^n_{\hat \tau} }  \left(  {\cal U}_{\hat \tau}   1_{  [0, \varepsilon]  }   \right) (x) \   \hat \mu_A (dx) 
  \le      
  \hat \gamma_2   \int_{W}  \left( 
         \int_{A}   
        p({\bm t}, x)  \, 1_{ [ 0,  \tau_2^{-n} (\varepsilon) ] } 
                         ( T_{{\bm t}} (x))    \, m(dx)  
         \right) \nu (d {\bm t})  
\end{equation}
for any small $ \varepsilon >0 $.
In a way similar to obtaining (7.3),  putting 
$
C_7 =  p_1^\ast \, C_2  / M_1
$,
we have
\begin{eqnarray}
 \int_{W}   \int_{A}   
        p({\bm t}, x)  \, 1_{ [ 0,  \tau_2^{-n} (\varepsilon) ] } 
                         ( T_{{\bm t}} (x))  \, m(dx)   \nu (d {\bm t}) 
& \leq  & C_7  \int_{ a_1 }^{ b_1 }  \left( \tau_2^{-n} (\varepsilon) \right)^s ds \\                     
& \leq  &
\frac{C_7}{ \log  \tau_2^{-n} (\varepsilon)  } 
  \left( \left( \tau_2^{-n} (\varepsilon) \right)^{ b_1 } -  \left( \tau_2^{-n} (\varepsilon) \right)^{ a_1} \right).  \nonumber 
\end{eqnarray}
It is easy to see that 
$$
C_{\tau_{2,1}} \left(  \varepsilon^{-(a_0 + \kappa -1)} + n  \right)^{- 1/(a_0  + \kappa -1)} 
\le 
 \tau_2^{-n} (\varepsilon) 
\le 
C_{\tau_{2,2}}  \left(  \varepsilon^{-(a_0 + \kappa -1)} + n  \right)^{- 1/(a_0  + \kappa -1)}  
$$
 for small $\varepsilon >0 $, where $ 0 < C_{\tau_{2,1}} <1 $ and $C_{\tau_{2,2}}  $ are constants.
So, we have 
\begin{eqnarray*}
 \int_{ a_1 }^{ b_1 }  \left( \tau_2^{-n} (\varepsilon) \right)^s ds 
& \le & 
\frac  {  C_{\tau_{2,2}}   \left(  \varepsilon^{-(a_0 + \kappa -1)} + n  \right)^{- a_1 / (a_0   + \kappa -1)} }
     {- \log C_{\tau_{2,1}} +  \frac 1 { a_0 + \kappa -1} \log  \left(  \varepsilon^{-(a_0 + \kappa -1)} + n  \right)}   \\
& \le & 
\frac  {  C_{\tau_{2,2}}   \left(  \varepsilon^{-(a_0 + \kappa -1)} + n  \right)^{- a_1 / (a_0   + \kappa -1)} }
     {  \frac 1 { a_0 + \kappa -1} \log  \left(  \varepsilon^{-(a_0 + \kappa -1)} + n  \right)}
\end{eqnarray*}
for any  small  $\varepsilon >0 $. 
Thus, using  $ a_1 / (a_0  + \kappa -1)  > 1$,  we have 
\begin{eqnarray} 
\sum_{n=0}^{\infty} 
  \int_{ a_1 }^{ b_1 }  \left( \tau_2^{-n} (\varepsilon) \right)^s ds 
& \le &
  \sum_{n=0}^{\infty}     
 \frac  { C_{\tau_{2,2}}    \left(  \varepsilon^{-(a_0 + \kappa -1)} + n  \right)^{  - a_1 / (a_0   + \kappa -1)} }
     { \frac 1 { a_0 + \kappa -1 } \log  \left(  \varepsilon^{-(a_0 + \kappa -1)} + n  \right)   }     \nonumber  \\ 
 & \le &
    \sum_{n=0}^{\infty}     
 \frac  {  C_{\tau_{2,2}}    \left(  \varepsilon^{-(a_0 + \kappa -1)} + n  \right)^{- a_1  / (a_0   + \kappa -1)} }
     { -  \log    \varepsilon  }     \nonumber   \\
 & \le &
  \frac {   2 C_{\tau_{2,2}} (a_0 + \kappa - 1 ) \varepsilon^{a_1 + 1 - a_0  - \kappa}} 
                  { - ( a_1 + 1 - a_0 - \kappa )  \log    \varepsilon   }.
\end{eqnarray}
Put 
$C_8 =  2 C_7 \, C_{\tau_{2,2}} (a_0 + \kappa - 1 ) /   (a_1 - a_0 - \kappa +1 )$.
 Then, by Theorem 6.5 and the inequalities (7.6)-(7.9),   we obtain 
$$ \displaystyle 
\mu([0, \varepsilon ]) 
\leq   \frac{\gamma_2 }{ \hat \gamma_1  \inf_{x \in [0, \hat c)} \hat p (x)} \hat{\mu}( [0, \varepsilon]  )
\leq 
\frac{\gamma_2 \hat \gamma_2 C_8   \,  \varepsilon^{ a_1 + 1 - a_0  - \kappa} }
     {-  \hat \gamma_1  \inf_{x \in [0, \hat c)} \hat p (x) \cdot \log    \varepsilon}
$$
for any small $   \varepsilon > 0$.
Hence, the upper estimate of (2) is obtained and proof is completed. 
  \hfill $\Box$

\vspace{4ex}

\noindent 
{\it
Proof of Theorem  7.2.}

\noindent 
[Lower estimate]

Since $\nu_1$ is the counting measure, 
we obtain 
\begin{eqnarray*}
  \int_{ A }
             { \widetilde{{\cal U}}^n_{\bar \tau} }  \left(  {\cal U}_{\bar \tau}   1_{  [0, \varepsilon]  }   \right) (x) \ \bar \mu_A (dx) 
\geq 
 \bar \gamma_1  C_5
\left( \tau_1^{-n} (\varepsilon) \right)^{a_1} 
\end{eqnarray*}
instead of the inequalities  (7.2) and (7.3).
Since 
$
 \tau_1^{-n} (\varepsilon) \ge C_{\tau_{1 ,1} }  \left(  \varepsilon^{-(a_0 -1)} + n  \right)^{- 1/(a_0 -1)}  
$ for small $\varepsilon >0 $, 
we have 
$$
  \left( \tau_1^{-n} (\varepsilon) \right)^{a_1}  
\ge 
C_{\tau_{1 ,1}}  \left(  \varepsilon^{-(a_0-1)} + n  \right)^{- a_1  /(a_0-1)}.
$$

(1) 
Assume 
$ a_0 \ge  a_1 +1$. 
Then, in a way similar to the proof of Theorem 7.1 (1), 
we obtain
$
\mu ([0, \varepsilon]) = \infty 
$ 
for any small $ \varepsilon > 0$. 

(2)$^\ast$ 
Assume $ a_0  <  a_1 +1$. Then, in a way similar to the proof of  the lower estimate of Theorem 7.1 (2), 
we obtain
$$ \displaystyle 
\bar \mu ([0, \varepsilon]) \ge 
 \bar  \gamma_1 C_9 \,  \varepsilon^{ a_1 + 1 - a_0  }
 \quad \mbox{for any small} \   \varepsilon > 0, 
$$
where 
$ \displaystyle 
C_9 = 
\frac { C_5 C_{\tau_{1 ,1}} (a_0 -1) } 
    { 2 (a_1+1  - a_0 )}
$.
Thus, it follows from Theorem 6.2 that 
$$
\mu ([0, \varepsilon]) 
\ge  \frac { \gamma_1 \bar \gamma_1 C_9  \varepsilon^{a_1 + 1 - a_0 } } { \bar \gamma_2}
\quad \mbox{for any small} \   \varepsilon > 0 
  \quad \mbox {if}  \quad    a_0  \le a_1 +1, 
$$
which shows the lower estimate of (2)$^\ast$. 


\vspace{1ex}

\noindent 
[Upper estimate]

Since $\nu_1$ is the counting measure, 
we obtain 
$$
 \int_{ A }
             { \widetilde{{\cal U}}^n_{\hat \tau} }  \left(  {\cal U}_{\hat \tau}   1_{  [0, \varepsilon]  }   \right) (x) \,  \hat \mu_A (dx)
\le 
2 \hat \gamma_2  C_7
\left( \tau_2^{-n} (\varepsilon) \right)^{ a_1} 
\quad \mbox{for any small} \ \varepsilon >0 
$$
instead of the inequalities (7.7) and (7.8).
By a minor modification of the proof of the upper estimate of Theorem 7.1, 
we have 
$$ \displaystyle 
\mu([0, \varepsilon ]) 
\leq   \frac{\gamma_2 }{ \hat \gamma_1  \inf_{x \in [0,c)} \hat p (x)} \hat{\mu}( [0, \varepsilon]  )
\leq 
\frac{ 2 \gamma_2 \hat \gamma_2 C_7   \,  \varepsilon^{a_1 +  1 - a_0  - \kappa} }
     {-  \hat \gamma_1  \inf_{x \in [0,c)} \hat p (x)  }
$$
for any small $   \varepsilon > 0$, 
which shows the upper  estimate of (2)$^\ast$. 
  \hfill $\Box$

\vspace{5ex}

Now, we are going to apply our comparison theorems to show Theorems 7.5 and 7.6. 

\vspace{3ex}

\noindent 
{\it 
Proof of Theorem 7.5.
}

\noindent 
[Lower estimate]

Let  $\bar  c$ be a small constant with $0 < \bar c < \min \{ c,   \varepsilon_0 \}$.
Let $\tau_1 :[0, c ) \to [0, 1) $ be a surjective and  $C^1$ map 
such that 
$\tau_1(x) = b_0 x  $ on $ [0,   \bar c ] $ and that  $ \tau_1^{\prime} (x) > 1$ on $ [0, c)$.
Then,
for each $ \bm t = (t, s) \in W$ 
$$  T_{ \bm t }  (0,  \varepsilon ) \subset \tau_1 (0,  \varepsilon ) \ \mbox{for any $ \varepsilon \in (0,  \bar c ) $}.
$$
For each $\bm t \in W$,
define $ \bar{\tau}_{\bm t} : [0, 1]  \to [0,1]$ 
by 
$$
\bar{\tau}_{\bm t} (x) := 
\left\{
\begin{array}{ll}  
\tau_1 (x),   &    \quad  [0, c),  \\   
T_{\bm t}(x), &  \quad  [ c, 1].
\end{array}
 \right.
$$
Then, by Theorem 4.1,  
the random map 
$ \bar{\tau} = \{ \bar {\tau}_{\bm t} , p(\bm t,x) : \bm t \in W \}$
has an absolutely continuous invariant measure $ \bar {\mu}$ with $\bar \mu([c,1]) =1 $.
Let $\bar h$ be a density of $ \bar {\mu}$.
By Theorem 4.1 and Lemmas 7.8, 7.9 and 7.10, we have constants $\bar \gamma_1 $ and $\bar \gamma_2$ such that 
$ 0<  \bar \gamma_1 \le \bar{h} (x) \le \bar \gamma_2 < \infty$ 
on $  \bigcup_{\bm t \in W} T_{\bm t} |_{[c,1]}^{-1}(B(0)) $.

As we obtained  (7.2) and (7.3) in the proof of  the lower estimate of Theorem 7.1, 
we obtain 
\begin{eqnarray*}
& &   \int_{ A }
             { \widetilde{{\cal U}}^n_{\bar \tau} }  \left(  {\cal U}_{\bar \tau}   1_{  [0, \varepsilon]  }   \right) (x) \ \bar \mu_A (dx)
 \geq  \bar \gamma_1 C_5 \int_{ a_1 }^{ b_1 }  \left( \tau_1^{-n} (\varepsilon) \right)^s ds
              \\
 &  = &   
\frac{ \bar \gamma_1 C_5}{ \log  \tau_1^{-n} (\varepsilon)  } 
  \left( \left( \tau_1^{-n} (\varepsilon) \right)^{ b_1 } -  \left( \tau_1^{-n} (\varepsilon) \right)^{ a_1 } \right).
\end{eqnarray*}
Especially, we have 
\begin{equation}
 \int_{ A }
             {\cal U}_{\bar \tau}   1_{  [0, \varepsilon]  }  (x) \  \bar \mu_A (dx)
             \ge 
             \frac{ \bar \gamma_1 C_5 (   \varepsilon^{ b_1 } -  \varepsilon^{ a_1 }  )}{ \log  \varepsilon   } 
             \ge 
             \frac{ \bar \gamma_1 C_5  \varepsilon^{ a_1 } }{ - 2  \log  \varepsilon }
             \quad \mbox{for any small} \  \varepsilon >0.
\end{equation}
Thus,  we obtain 
\begin{eqnarray*}
\sum_{n=0}^{\infty} 
         \int_{ A }     { \widetilde{{\cal U}}^n_{\bar \tau} }  \left(  {\cal U}_{\bar \tau}   1_{  [0, \varepsilon]  }   \right) (x) \   \bar \mu_A (dx)
& \ge & 
\int_{ A }
             {\cal U}_{\bar \tau}   1_{  [0, \varepsilon]  }  (x) \   \bar \mu_A (dx)\\
& \ge &     
 \frac{ \bar \gamma_1  C_5  \, \varepsilon^{ a_1 } }{ - 2  \log  \varepsilon }, 
\end{eqnarray*}
which means 
$$
\bar \mu ([0, \varepsilon]) 
\ge    
\frac  { \bar \gamma_1 C_5 \,  \varepsilon^{ a_1} }{ -  2 \log \varepsilon} \quad \mbox{for any small} \  \varepsilon >0.          
$$
Therefore, by Theorem 6.2 we obtain 
\begin{equation}
\mu ([0, \varepsilon])
\ge 
\frac  { \gamma_1 \bar \gamma_1  C_5 \, \varepsilon^{ a_1} }{ -  2  \bar \gamma_2 \log \varepsilon}  
\quad \mbox{for any small} \  \varepsilon >0, 
\end{equation}
which shows the lower estimate in the theorem.

\vspace{1ex}

\noindent 
[Upper estimate]

Let  $\hat  c$ be a small constant with $0 < \hat c < \min \{ c,   \varepsilon_0 \}$. 
Let $\tau_2 :[0, c ) \to [0, 1) $ be a surjective and  $C^1$ map 
such that 
 $\tau_2(x) = a_0 x $ on $ [0,   \hat c ] $ and that  $ \tau_2^{\prime} (x) > 1$ on $ [0, c)$.
Put 
$$ W_2 (\tau_2)= \{ \bm t \in W \ 
 | \  T_{ \bm t }  (0,  \varepsilon ) \supset \tau_2 (0,  \varepsilon ) \ \mbox{for any $0 < \varepsilon < \hat c$}
\}
$$
 and 
$$
\hat p (x) :=
\int_{W_2(\tau_2)}  p( \bm t, x) \nu(d \bm t) \quad \mbox{for $ x \in  [0,  c) $}. 
$$
Then 
$
\hat p (x) = 1 
$
for $ x \in  [0,  c) $.
For each $ \bm t \in W$,
define $ \hat{\tau}_{\bm t} : [0, 1]  \to [0,1]$ 
by 
$$
\hat{\tau}_{\bm t} (x) := 
\left\{
\begin{array}{ll}  
\tau_2 (x),   &    \quad  [0, c),  \\   
T_{\bm t} (x), &  \quad  [ c, 1].
\end{array}
 \right.
$$
By Theorem 4.1,  the random map 
$ \hat{\tau} = \{ \hat{\tau}_{\bm t} , p( \bm t, x) ;  \bm t \in W \}$
has an absolutely continuous invariant measure $ \hat{\mu}$
with $ \hat \mu([c,1]) =1 $.
Let $\hat h$ be a density of $ \hat{\mu}$.
As $ \bar h$ in the lower estimate, it follows  that  $ 0< \hat \gamma_1 \le  \hat{h}(x)  \le \hat \gamma_2 < \infty$
on $  \bigcup_{\bm t \in W} T_{\bm t } |_{[c,1]}^{-1}(B(0)) $, 
 where $\hat \gamma_1 $ and $\hat \gamma_2 $ are constants.

As we obtained (7.7) and (7.8) 
in the proof of the upper estimate of Theorem 7.1, 
we obtain 
\begin{eqnarray}
& &   \int_{ A }
             { \widetilde{{\cal U}}^n_{\hat \tau} }  \left(  {\cal U}_{\hat \tau}   1_{  [0, \varepsilon]  }   \right) (x) \  \hat \mu_A (dx)
 \leq  \hat  \gamma_2 C_7  \int_{a_1 }^{ b_1 }  \left( \tau_2^{-n} (\varepsilon) \right)^s ds
              \\
 &  = &   
\frac{ \hat  \gamma_2  C_7}{ \log  \tau_2^{-n}  (\varepsilon)  } 
  \left( \left( \tau_2^{-n} (\varepsilon) \right)^{b_1} -  \left( \tau_2^{-n} (\varepsilon) \right)^{ a_1 } \right).  \nonumber 
\end{eqnarray}
Since 
$
 \tau_2^{-n} (\varepsilon) =  a_0^{-n} \varepsilon
$, 
we have 
\begin{equation}
\int_{a_1 }^{ b_1}  \left( \tau_1^{-n} (\varepsilon) \right)^s ds
=
\frac  { a_0^{- a_1 n } \varepsilon^{ a_1 } - a_0^{- b_1 n } \varepsilon^{ b_1 }}{ n \log a_0 - \log \varepsilon}
\leq
\frac  { a_0^{-a_1  n } \varepsilon^{ a_1} }{ - \log \varepsilon}. 
\end{equation}
Thus, by (7.12) and (7.13) 
we have
$$
\hat \mu ([0, \varepsilon]) = 
\sum_{n=0}^{\infty} 
  \int_{ A }
             { \widetilde{{\cal U}}^n_{\hat \tau} }  \left(  {\cal U}_{\hat \tau}   1_{  [0, \varepsilon]  }   \right) (x) \  \hat \mu_A (dx) 
\le
\frac  { \hat  \gamma_2  C_7 \,  \varepsilon^{ a_1} }{ - ( 1 - a_0^{- a_1} ) \log \varepsilon}.             
$$
Therefore, remarking that  $ \hat p (x) = 1$,   it follows from Theorem 6.5 that 
\begin{equation}
 \mu ([0, \varepsilon]) \le
 \frac {\gamma_2}{ \hat \gamma_1 }\hat \mu ([0, \varepsilon]) 
\le
\frac  { \gamma_2 \hat  \gamma_2  C_7 \,  \varepsilon^{ a_1} }
        { -  \hat \gamma_1 ( 1 - a_0^{- a_1} ) \log \varepsilon}
    \quad \mbox{for any small $ \varepsilon > 0$}.
\end{equation}
Hence, the upper  estimate in the theorem is obtained.
  \hfill $\Box$

\vspace{3ex}

{\it 
\noindent 
Proof of Theorem 7.6.}

\noindent 
[Lower estimate]

Similarly to obtaining (7.10) in the proof of Theorem 7.5, we obtain 
$$
 \int_{ A }
             {\cal U}_{\bar \tau}   1_{[0, \varepsilon ] }  (x) 
              \bar \mu_A (dx)
             \ge \bar \gamma_1  C_5 \varepsilon^{a_1}
              \quad \mbox{for any small} \  \varepsilon >0.
$$
Using this inequality  instead of (7.10),  similarly to obtaining (7.11) in the proof of Theorem 7.5,  we obtain 
$$
\mu ([0, \varepsilon])
\ge 
\frac  { \gamma_1 \bar \gamma_1  C_5 \, \varepsilon^{ a_1} }{  \bar \gamma_2 }  
\quad \mbox{for any small} \  \varepsilon >0. 
$$
Hence, the lower estimate in the theorem is obtained. 

\vspace{1ex}

\noindent 
[Upper estimate]

Similarly to  obtaining (7.12) and (7.13) in the proof of  the upper estimate of Theorem 7.5, we obtain 
$$
  \int_{ A }
             { \widetilde{{\cal U}}^n_{\hat \tau} }  \left(  {\cal U}_{\hat \tau}   1_{  [0, \varepsilon]  }   \right) (x) \   \hat \mu_A (dx)
\le 
2  \hat  \gamma_2 C_7 
  \left( \tau_2^{-n} (\varepsilon) \right)^{ a_1 }
 =  2  \hat  \gamma_2 C_7 a_0^{-a_1 n }  \varepsilon^{ a_1 }
   \quad \mbox{for any small} \  \varepsilon >0.
$$
Using the above inequality  instead of (7.12) and (7.13),  similarly to obtaining (7.14) in the proof of Theorem 7.5,  we obtain 
$$
 \mu ([0, \varepsilon]) \le
 \frac {\gamma_2}{ \hat \gamma_1 }\hat \mu ([0, \varepsilon]) 
\le
\frac  { 2 \gamma_2 \hat  \gamma_2  C_7 \,  \varepsilon^{ a_1} }
        {  \hat \gamma_1 ( 1 - a_0^{- a_1} ) }    \quad \mbox{for any small} \  \varepsilon >0.
$$
Hence, the upper estimate in the theorem is obtained. 
  \hfill $\Box$

\vspace{2ex}

 \section{Concluding remark} 
  
 Bisides the random maps we have just seen in the previous section,  
 there are a lot of  random maps to which our comparison theorems can  
  be applied. 
  
Even if we change some parts of the setting in the previous section, 
 we expect we can 
 apply  our comparison theorems and estimate  the absolutely continuous invariant measures of the random maps.
  
 Some  bigger changes are also possible. 
 For example,   
 we can also apply our comparison theorems to random piecewise convex  maps, which is discussed in [InT], 
 where we show the existence of absolutely continuous invariant measures of  random piecewise convex  maps, 
 and apply  our comparison theorems 
 to the estimate of the  invariant measures.
On deterministic systems,
the existence of absolutely continuous invariant measures of  piecewise convex maps were studied by Lasota and Yorke (see [LY2]),  
and the results of Lasota and Yorke were  generalized  in [In1] and [In2]. 
 
 We can  expect that  there are a lot of applications of our comparison theorems.

  \vspace{5ex}
  
 \noindent   
{\bf Acknowledgments. }
 
 \noindent 
The author would like to express his sincere gratitude to Dr. H. Toyokawa (Kitami Institute of Technology)
who gave the author  valuable questions 
and with whom he had  fruitful  discussions.
 
\vspace{10ex}

\begin{center}
{\bf REFERENCES}
\end{center}

{\footnotesize 

\newcommand{\namelistlabel}[1]{\mbox{#1}\hfil}
\newenvironment{namelist}[1]{%
\begin{list}{}
  {\let\makelabel\namelistlabel
   \settowidth{\labelwidth}{#1}
   \setlength{\leftmargin}{1.1\labelwidth}}
}{%
\end{list}}

\begin{namelist}{xxxxxxxx}

\item[{[BaBD]}]  W. Bahsoun, C. Bose and Y. Duan,
{\em Decay  of correlation for random intermittent maps}, 
{Nonlinearity}  {27} (2014), 1543-1554.

\item[{[BaG]}] 
W. Bahsoun and P. G\'ora, 
{\em Position dependent random maps in one and higher dimensions},
{ Studia Math.} {166} (2005), 271-286.

\item[{[BoG]}]  A. Boyarsky and P. G\'ora,
{\em Laws of Chaos},  Birkh\"auser, Boston, 1997. 

\item[{[CHMV]}]
G.  Cristadoro, N.  Haydn, P.  Marie, 
and S. Vaienti, 
{\em Statistical properties of intermittent maps with
unbounded derivative},  {Nonlinearity}  {23} (2010), 1071-1095.

\item[{[GBo]}]  
 P. G\'ora and  A. Boyarsky,
{\em Absolutely continuous invariant measures for
random maps with position dependent probabilities}, 
{  J. Math. Anal. Appl.} {278}  (2003), 225-242.

\item[{[In1] }]  T. Inoue,  {\em Asymptotic stability of densities for piecewise convex maps,}
{Ann. Polon. Math. } {LVII} (1992), 83-90. 

\item[{[In2] }]  T. Inoue,  {\em Weakly attracting repellors for piecewise convex maps,}
{Japan J. Indust. Appl. Math.} {9} (1992), 413-430.

\item[{[In3] }]  T. Inoue,  
{\em Invariant measures for position dependent random maps 
with continuous random parameters},  
Studia Math. 208  (2012), 11-29.

\item[{[In4]}]  
T. Inoue, 
{\em First return maps of random maps 
and invariant measures,}
Nonlinearity {33}  (2020),
249-275.

\item[{[InT]}]  
T. Inoue. and H. Toyokawa, 
{\em Invariant measures for random piecewise convex  maps}
Preprint.

\item[{[KiLi]}]
Y. Kifer and P. D. Liu, 
{\em Random dynamics },
Handbook of dynamical systems, {1B,}  379-499, 
Elsevier B. V., Amsterdam, 
2006.

\item[{[Ko]}] 
Z.S. Kowalski, 
{\em Invariant measure  for piecewise monotonic transformation has a positive lower bound on its support }
Bull. Acad. Polon. Sci. Ser. Sci. Math. {27} (1979), 53-57.

\item[{[LM]}]  A. Lasota and M. C. Mackey, {\em Chaos, fractals, and noise,  
Stochastic aspects of dynamics},
 Appl. Math. Sci.  { 97},  Springer, New York, Berlin, Heidelberg, 1994.  

\item[{[LY1]}]  A. Lasota and J. A. Yorke,
{\em On the existence of Invariant densities for piecewise monotonic transformations},
Trans. Amer. Math. Soc. {186}
(1973), 481-488.

\item[{[LY2]}] 
A. Lasota and J. A. Yorke, 
{\em Exact dynamical systems and the Frobenius-Perron operator}, Trans. Amer. Math. Soc., 273 (1982),  375-384.

\item[{[LSV]}] C. Liverani, B. Saussol and S. Vaienti, 
{\em A probabilistic approach to intermittency},  
Ergodic Theory Dynam. Systems {19} (1999), 671-685.

\item[{[MaSU]}]   V. Mayer, B. Skorulski, and M. Urbanski, 
{\em Distance expanding random mappings, thermodynamical formalism, Gibs measures and 
 fractal geometry},   Lecture Notes in Math.  2036. Springer, Heidelberg, 2011.

\item[{[Mo]}]  T. Morita,   {\em Random iteration of one-dimensional transformations}, 
Osaka J. Math. {22} (1985), 489-518.

\item[{[Pe]}]  S. Pelikan, {\em Invariant densities for random maps of the  interval},
Trans. Amer. Math. Soc. {281}
(1984), 813-825.

\item[{[Pi]}] G. Pianigiani, {\em First return map and invariant measures}, 
Israel J. Math. {35} (1980),  32-48.

\item[{[R]}]  M. Rychlik,  {\em Bounded variation and invariant measures}, {Studia Math.} {76}  (1983), 69-80.

\item[{[T1]}]  M. Thaler,  {\em Estimates of the invariant densities of endomorphisms 
with indifferent fixed points}, 
Israel J.Math. {37} (1980),  304-314. 

\item[{[T2]}]  M. Thaler,  {\em Transformations on $[0,1]$ with infinite invariant measures}, 
Israel J.Math. {46} (1983),  67-96.

\end{namelist}

}

\end{document}